\documentclass[12pt]{amsart}
\usepackage{amssymb}
\usepackage{amsmath}
\usepackage{amscd}
\usepackage{graphicx}
\theoremstyle{plain}
\newtheorem{theorem}{Theorem}[section]
\newtheorem{corollary}[theorem]{Corollary}
\newtheorem{proposition}[theorem]{Proposition}
\newtheorem{lemma}[theorem]{Lemma}
\newtheorem{problem}[theorem]{Problem}
{\theoremstyle{remark}

\newtheorem{remark}[theorem]{Remark}}
{\theoremstyle{definition}
\newtheorem{definition}[theorem]{Definition}

\newtheorem{example}[theorem]{Example}}

\newcommand{\rom}{\renewcommand{\labelenumi}{{\rm (\roman{enumi})}}%
\renewcommand{\itemsep}{0pt}}

\newcommand{\N}{\mathbb{N}}
\newcommand{\Z}{\mathbb{Z}}
\newcommand{\R}{\mathbb{R}}
\newcommand{\C}{\mathbb{C}}
\newcommand{\T}{\mathbb{T}}
\newcommand{\Q}{\mathbb Q}
\newcommand{\e}{\varepsilon}
\newcommand{\cK}{{\mathcal K}}
\newcommand{\cL}{{\mathcal L}}

\newcommand{\M}{\mathbb{M}}

\newcommand{\cO}{{\mathcal O}}
\newcommand{\cT}{{\mathcal T}}
\newcommand{\ip}[2]{\langle{#1},{#2}\rangle}

\newcommand{\s}[3]{{{#1}^{#2}_{\textnormal{#3}}}}
\newcommand{\mat}[4]{\bigg(\begin{array}{cc}
{#1}&{#2}\\{#3}&{#4}\end{array}\bigg)}
\newcommand{\rs}[1]{{\rm{\textnormal #1}}}
\newcommand{\td}{{\tilde d}}
\newcommand{\tr}{{\tilde r}}
\newcommand{\tE}{\widetilde{E}}

\newcommand{\Ca}{$C^*$-al\-ge\-bra }
\newcommand{\CA}{$C^*$-al\-ge\-bra}
\newcommand{\Cas}{$C^*$-al\-ge\-bras }
\newcommand{\shom}{$*$-ho\-mo\-mor\-phism }

\newcommand{\Cc}{$C^*$-cor\-re\-spon\-dence }
\newcommand{\Ccs}{$C^*$-cor\-re\-spon\-dences }
\newcommand{\Csa}{$C^*$-sub\-al\-ge\-bra }

\DeclareMathOperator{\id}{id}

\DeclareMathOperator{\coker}{coker}
\DeclareMathOperator{\Orb}{Orb}

\DeclareMathOperator{\cspa}{\overline{span}}

\begin{document}
\title[A class of $C^*$-algebras IV, pure infiniteness]
{A class of \boldmath{$C^*$}-algebras generalizing 
both graph algebras and homeomorphism \boldmath{$C^*$}-algebras IV, \\
pure infiniteness.}
\author[Takeshi KATSURA]{Takeshi KATSURA}
\address{Department of Mathematics, 
Hokkaido University, Kita 10, Nishi 8, 
Kita-Ku, Sapporo, 060-0810, JAPAN}
\email{katsura@math.sci.hokudai.ac.jp}
\date{}


\maketitle

\begin{abstract}
This is the final one in the series of papers 
where we introduce and study 
the \CA s associated with topological graphs. 
In this paper, 
we get a sufficient condition on topological graphs 
so that the associated $C^*$-algebras are simple and 
purely infinite. 
Using this result, 
we give one method to construct all Kirchberg algebras 
as $C^*$-algebras associated with topological graphs. 
\end{abstract}

\setcounter{section}{-1}

\section{Introduction}\label{Intro}

The classification theory of simple separable nuclear \CA s 
by $K$-theory has been rapidly developed recently. 
This classification was completed for the purely infinite case 
independently by Kirchberg \cite{Ki} and Phillips \cite{P}. 
Recall that a simple \Ca is said to be {\em purely infinite} 
if every non-zero hereditary \Csa has an infinite projection. 

\medskip
\noindent 
{\bf Definition}.
A {\em Kirchberg algebra} is a simple, separable, nuclear, purely infinite 
\Ca satisfying the universal coefficient theorem of \cite{RoSc}. 
\medskip

\noindent 
{\bf Theorem} (Kirchberg, Phillips). 
Two non-unital Kirchberg algebras $A$ and $B$ are isomorphic 
if and only if $(K_0(A),K_1(A))\cong (K_0(B),K_1(B))$, 
and two unital Kirchberg algebras $A$ and $B$ are isomorphic 
if and only if $(K_0(A),[1_A],K_1(A))\cong (K_0(B),[1_B],K_1(B))$. 
\medskip

See R\o rdam's book \cite{RS} for detailed definitions 
and a proof of the above theorem 
(note that in \cite{RS} a Kirchberg algebra is not 
assumed to satisfy the universal coefficient theorem). 
There are many ways to construct Kirchberg algebras, 
and all pairs $(G_0,G_1)$ of countable abelian groups 
appear as $K$-groups of non-unital Kirchberg algebras, 
and for all $g\in G_0$ 
there exists a unital Kirchberg algebra $A$ 
with $(K_0(A),[1_A],K_1(A))\cong (G_0,g,G_1)$ 
(see \cite[Subsection 4.3]{RS}). 

One way to construct Kirchberg algebras 
is by graph algebras. 
A {\em graph algebra} is a \Ca associated with a (directed) graph, 
which generalizes Cuntz-Krieger algebras defined in \cite{CK} 
(see \cite{Ra} for a definition and properties of graph algebras). 
We know a necessary and sufficient condition on graphs 
so that the associated graph algebras are simple and purely infinite, 
and in this case the graph algebras become Kirchberg algebras 
(see \cite[Remark 4.3]{Ra}).  
However we cannot construct all Kirchberg algebras 
by graph algebras 
because $K_1$-groups of graph algebras are always free. 
Spielberg constructed all non-unital Kirchberg algebras 
mixing the constructions of graph algebras and 
the higher rank graph algebras \cite{Sp}. 

In this paper, 
we examine for which topological graphs 
the associated \CA s are simple and purely infinite. 
We also construct all Kirchberg algebras 
as \CA s associated with topological graphs. 
A {\em topological graph} is a quadruple $E=(E^0,E^1,d,r)$ 
consisting of the set of vertices $E^0$ 
and the set of edges $E^1$ which are locally compact spaces, 
and two continuous maps $d,r\colon E^1\to E^0$ 
which indicate the domains and the ranges of edges, respectively. 
The domain map $d$ is assumed to be locally homeomorphic. 
When $E^0$ is discrete, 
a topological graph $E=(E^0,E^1,d,r)$ 
is nothing but an ordinary graph, 
which we call a {\em discrete graph} in this paper. 
In \cite{Ka4}, 
we introduce a way to associate a \Ca $\cO(E)$ 
to a topological graph $E$, 
which uses a construction of \CA s from \Ccs 
(see \cite{Ka7}). 
When a topological graph $E$ is a discrete graph, 
the \Ca $\cO(E)$ is isomorphic to the graph algebra of $E$. 

The author thinks that topological graphs are 
kinds of dynamical systems. 
The triple $(E^1,d,r)$ can be considered as a continuous 
multi-valued map from $E^0$ to itself 
(which we call a topological correspondence in \cite{Ka4}). 
For each natural number $n\in\N=\{0,1,2,\ldots\}$, 
the triple $(E^n,d^n,r^n)$ consisting of 
the set $E^n$ of paths with length $n$ 
and the domain and range maps $d^n,r^n\colon E^n\to E^0$ 
is the $n$-times iteration of $(E^1,d,r)$, 
and the map $\N\ni n\mapsto (E^n,d^n,r^n)$ 
defines a (multi-valued) action of the semigroup $\N$ on $E^0$. 
We can think that the \Ca $\cO(E)$ is a crossed product of this action. 
There has been many works on constructing Kirchberg algebras 
by crossed products (for example \cite{A}, \cite{LS}), 
but they used non-amenable groups 
because (ordinary) crossed products of commutative \CA s by 
amenable groups never be Kirchberg algebras. 
Our work in this paper is related to 
these studies 
although the \CA s associated with topological graphs 
are regarded as crossed products by 
the abelian semigroup $\N$. 
We note that 
singly generated dynamical systems 
introduced by Renault in \cite{Re} 
are examples of topological graphs, 
and the construction of \CA s in \cite{Re} 
from singly generated dynamical systems using groupoids 
is the same as the one here 
considering them as topological graphs. 
For these \CA s, 
there were some works on pure infiniteness 
(see for example \cite{Hj}), 
and it seems that our work here has a large overlap 
with them (cf.\ \cite{Ka8}). 

We are going to explain the two main results in this paper. 
In \cite{Ka6}, 
we extended many notions of dynamical systems 
such as orbits and minimality to topological graphs, 
and got equivalent conditions on topological graphs $E$ 
for the associated \CA s $\cO(E)$ to be simple 
(see Proposition \ref{PropSimple}). 
In Section \ref{SecCont} of this paper, 
we introduce the notion of contracting topological graphs 
(Definition \ref{Def(PI)}), 
and prove the following theorem which is 
the first main theorem of this paper. 

\medskip
\noindent 
{\bf Theorem A}. 
For a minimal and contracting topological graph $E$, 
the \Ca $\cO(E)$ is simple and purely infinite. 
\medskip

A topological graph $E=(E^0,E^1,d,r)$ 
is said to be {\em second countable} 
if both $E^0$ and $E^1$ are second countable. 
A topological graph $E$ is second countable 
if and only if $\cO(E)$ is separable 
(\cite[Proposition 6.3]{Ka4}), 
and in this case $\cO(E)$ satisfies the universal coefficient theorem 
(\cite[Proposition 6.6]{Ka4}). 
Since a \Ca associated with a topological graph 
is always nuclear (\cite[Proposition 6.1]{Ka4}), 
we get the following corollary of Theorem A. 

\medskip
\noindent 
{\bf Corollary B}.
For a second countable, minimal, contracting topological graph $E$, 
the \Ca $\cO(E)$ is a Kirchberg algebra. 
\medskip

The other main theorem of this paper is the following 
which follows from Proposition \ref{PropSt} 
and Proposition \ref{PropUn}. 

\medskip
\noindent 
{\bf Theorem C}. 
All Kirchberg algebras appear 
as \CA s of topological graphs. 
\medskip

Although there had been already many ways 
to construct Kirchberg algebras, 
it is important to give a new construction 
in order to attack some open problems. 
For example, 
construction of Kirchberg algebras as graph algebras 
was used to show that many Kirchberg algebras are 
semiprojective (see \cite{Sp2,Sz3}). 
In Appendix \ref{app3}, 
we discuss the possibility that our construction of Kirchberg algebras 
may extend the known results on the semiprojectivity of Kirchberg algebras 
(Remark \ref{remsemiproj}). 
For another example, 
a new construction of Kirchberg algebras 
may help to produce actions on them. 
In \cite{Sp3}, 
Spielberg used the construction in \cite{Sp} 
to show that all prime-order automorphisms of 
$K$-groups of Kirchberg algebras 
are induced by automorphisms of Kirchberg algebras 
having the same order. 
In \cite{Ka9}, 
we will extend his result using our construction. 

This paper is organized as follows. 
In Section \ref{SecPre}, 
we recall definitions of topological graphs 
and the \CA s associated with them. 
We also recall the 6-term exact sequence of $K$-groups 
of such \CA s 
and the criterion for their simplicity 
proved in \cite{Ka6}. 
In Section \ref{SecCont}, 
we introduce the notion of contracting topological graphs, 
and prove Theorem A. 
Section \ref{SecRem} is devoted 
to give some examples of contracting topological graphs 
and some remarks. 
In Section \ref{SecEnm} and Section \ref{SecEnmT}, 
we give a method to create topological graphs 
such that the associated \CA s are Kirchberg algebras 
and their $K$-groups are computable. 
We note that a similar construction can be found in \cite{D}. 
In Section \ref{SecKirch}, 
we use the method in Section \ref{SecEnmT} 
to get all Kirchberg algebras 
as \CA s associated with topological graphs, 
and thus prove Theorem C.

\section{Preliminaries}\label{SecPre}

We recall the definitions of topological graphs 
and their \CA s. 
For the detail, see \cite{Ka4}. 

\begin{definition}
A {\em topological graph} $E=(E^0,E^1,d,r)$ consists of 
two locally compact spaces $E^0$ and $E^1$, 
and two maps $d,r\colon E^1\to E^0$, 
where $d$ is locally homeomorphic 
and $r$ is continuous.
\end{definition}

From a topological graph $E=(E^0,E^1,d,r)$, 
we can define a \Cc $C_d(E^1)$ over $C_0(E^0)$ by 
\[
C_d(E^1)
=\{\xi\in C(E^1)\mid \ip{\xi}{\xi}\in C_0(E^0)\}, 
\]
where the inner product $\ip{\cdot}{\cdot}$ is defined by 
\[
\ip{\xi}{\eta}(v)=\sum_{e\in d^{-1}(v)}\overline{\xi(e)}\eta(e)
\]
for $\xi,\eta\in C_d(E^1)$ and $v\in E^0$, 
and the left and right actions are defined by 
\[
(f \xi g)(e)=f(r(e))\xi(e)g(d(e)) 
\]
for $f,g\in C_0(E^0)$, $\xi\in C_d(E^1)$ and $e\in E^1$. 
Note that the set of compact supported continuous functions $C_c(E^1)$ 
is contained in $C_d(E^1)$. 
The left action defines the \shom $\pi_r$ 
from $C_0(E^0)$ to the \Ca $\cL(C_d(E^1))$ 
of all adjointable operators on the Hilbert $C_0(E^0)$-module $C_d(E^1)$ 
such that $\pi_r(f)\xi=f\xi$. 

\begin{definition}
Let $E=(E^0,E^1,d,r)$ be a topological graph. 
We define an open subset $\s{E}{0}{rg}$ of $E^0$ by 
\begin{align*}
\s{E}{0}{rg}=\{v\in E^0\mid\ &
\text{there exists a neighborhood $V$ of $v$}\\
&\text{such that $r^{-1}(V)\subset E^1$ is compact, and $r(r^{-1}(V))=V$}\},
\end{align*}
and a closed subset $\s{E}{0}{sg}$ of $E^0$ by 
$\s{E}{0}{sg}=E^0\setminus \s{E}{0}{rg}$. 
\end{definition}

When $E^0$ is discrete, 
we have 
$\s{E}{0}{rg}=\{v\in E^0\mid 0<|r^{-1}(v)|<\infty\}$. 
The restriction of the \shom $\pi_r$ to the ideal $C_0(\s{E}{0}{rg})$ 
is an injection into the ideal $\cK(C_d(E^1))$ of $\cL(C_d(E^1))$ 
which is the closed span of the operators $\theta_{\xi,\eta}\in \cL(C_d(E^1))$ 
for $\xi,\eta\in C_d(E^1)$ 
defined by $\theta_{\xi,\eta}(\zeta)=\xi\ip{\eta}{\zeta}$. 

\begin{definition}\label{DefO(E)}
For a topological graph $E$, 
the $C^*$-algebra $\cO(E)$ is the universal \Ca 
generated by the images of a \shom $t^0\colon C_0(E^0)\to\cO(E)$ 
and a linear map $t^1\colon C_d(E^1)\to\cO(E)$ 
satisfying 
\begin{enumerate}
\rom
\item $t^1(\xi)^*t^1(\eta)=t^0(\ip{\xi}{\eta})$ for $\xi,\eta\in C_d(E^1)$, 
\item $t^0(f)t^1(\xi)=t^1(\pi_r(f)\xi)$ 
for $f\in C_0(E^0)$ and $\xi\in C_d(E^1)$,
\item $t^0(f)=\varphi(\pi_r(f))$ for $f\in C_0(\s{E}{0}{rg})$, 
\end{enumerate}
where $\varphi\colon \cK(C_d(E^1))\to \cO(E)$ is the \shom 
defined by $\varphi(\theta_{\xi,\eta})=t^1(\xi)t^1(\eta)^*$ 
for $\xi,\eta\in C_d(E^1)$, 
which is well-defined by (i). 
\end{definition}

The following proposition helps computations of the $K$-groups of $\cO(E)$. 

\begin{proposition}[{\cite[Corollary 6.10]{Ka4}}]\label{Kgroup}
For a topological graph $E$, 
the following sequence is exact; 
\[
\begin{CD}
K_0(C_0(\s{E}{0}{rg})) @>>\iota_*-[\pi_r]> K_0(C_0(E^0)) 
@>>t^0_*>  K_0(\cO(E)) \\
@AAA @. @VVV \\
K_1(\cO(E)) @<t^0_*<< K_1(C_0(E^0)) 
@<\iota_*-[\pi_r]<< K_1(C_0(\s{E}{0}{rg})),
\end{CD}
\]
where $\iota\colon C_0(\s{E}{0}{rg})\to C_0(E^0)$ is the inclusion map, 
and $[\pi_r]\colon K_i(C_0(\s{E}{0}{rg}))\to K_i(C_0(E^0))$ 
is the composition of the map 
$(\pi_r)_*\colon K_i(C_0(\s{E}{0}{rg}))\to K_i\big(\cK(C_d(E^1))\big)$ 
induced by the \shom $\pi_r\colon C_0(\s{E}{0}{rg})\to \cK(C_d(E^1))$ 
and the map $K_i\big(\cK(C_d(E^1))\big)\to K_i(C_0(E^0))$ 
induced by the Hilbert $C_0(E^0)$-module $C_d(E^1)$. 
\end{proposition}

We set $d^0=r^0=\id_{E^0}$ and $d^1=d, r^1=r$.
For $n=2,3,\ldots$, 
we recursively define a space $E^n$ of paths with length $n$ 
and domain and range maps $d^n,r^n\colon E^n\to E^0$ by
\[
E^n=\{(e,\mu)\in E^1\times E^{n-1}\mid d^1(e)=r^{n-1}(\mu)\},
\]
$d^n((e,\mu))=d^{n-1}(\mu)$ and $r^n((e,\mu))=r^1(e)$. 
For each $n\in\N$, 
$d^n$ is a local homeomorphism from $E^n$ to $E^0$ 
and $r^n$ is continuous. 
Note that there exists a natural isomorphism 
\[
E^{n+m}=\{(\mu,\nu)\in E^n\times E^m\mid d^n(\mu)=r^m(\nu)\}
\]
for $n,m\in\N$. 
For each $n\in\N$, 
we can define a \Cc $C_{d^n}(E^n)$ over $C_0(E^0)$ 
similarly as $C_d(E^1)$. 
This \Cc $C_{d^n}(E^n)$ is naturally isomorphic to 
the $n$-times tensor products of the \Cc $C_d(E^1)$, 
and using this fact we can define 
a linear map $t^n\colon C_{d^n}(E^n)\to \cO(E)$ 
such that 
\[
t^n(\xi)^*t^n(\eta)=t^0(\ip{\xi}{\eta}), \quad 
t^0(f)t^n(\xi)=t^n(\pi_{r^n}(f)\xi) 
\]
for $f\in C_0(E^0)$ and $\xi,\eta\in C_{d^n}(E^n)$ 
(see \cite[Section 2]{Ka4}). 
Recall that the norm of $\xi\in C_{d^n}(E^n)$ 
is defined by $\|\xi\|=\|\ip{\xi}{\xi}\|^{1/2}$, 
and the map $t^n$ is isometric. 

We recall from \cite{Ka6} conditions on a topological graph $E$ 
such that the \Ca $\cO(E)$ is simple. 
We set $E^*=\coprod_{n=0}^\infty E^n$. 
By extending $d^n$ and $r^n$, 
we get a local homeomorphism $d^*\colon E^*\to E^0$ 
and a continuous map $r^*\colon E^*\to E^0$.

\begin{definition}
We define the {\em positive orbit space} $\Orb^+(v)$ of $v\in E^0$ by 
\[
\Orb^+(v)=\{r^*(\mu)\in E^0\mid \mu\in E^*, d^*(\mu)=v\}.
\]
\end{definition}

An {\em infinite path} is a sequence $\mu=(e_1,e_2,\ldots,e_n,\ldots)$ 
with $e_i\in E^1$ and $d(e_i)=r(e_{i+1})$ for each $i=1,2,\ldots$. 
The set of all infinite paths is denoted by $E^\infty$. 
For an infinite path $\mu=(e_1,e_2,\ldots,e_n,\ldots)\in E^\infty$, 
we define its range $r^\infty(\mu)\in E^0$ to be $r(e_1)$. 

\begin{definition}
For $v\in E^0$, 
a {\em negative orbit} of $v$ 
is either a finite path $\mu\in E^n$ 
with $r^n(\mu)=v$ and $d^n(\mu)\in\s{E}{0}{sg}$, 
or an infinite path $\mu\in E^\infty$ with $r^\infty(\mu)=v$. 

For each negative orbit $\mu=(e_1,e_2,\ldots,e_n)\in E^n$ of $v\in E^0$ 
with $n\in\N\cup\{\infty\}$, 
the {\em negative orbit space} $\Orb^-(v,\mu)$ 
is defined by 
\[
\Orb^-(v,\mu)=\{v,d(e_1),d(e_2),\ldots,d(e_n)\}\subset E^0.
\]
\end{definition}

\begin{definition}
We define the {\em orbit space} $\Orb(v,e)$ of $v\in E^0$ 
with respect to a negative orbit $e$ of $v$ by 
\[
\Orb(v,\mu)=\bigcup_{v'\in\Orb^-(v,\mu)}\Orb^+(v').
\]
\end{definition}

\begin{definition}[{cf.\ \cite[Proposition 8.9]{Ka6}}]
A topological graph $E$ is said to be {\em minimal} 
if the orbit space $\Orb(v,\mu)$ is dense in $E^0$ 
for every $v\in E^0$ and every negative orbit $\mu$ of $v$. 
\end{definition}

A path $l=(e_1,\ldots,e_n)\in E^n$ for $n\geq 1$ is called 
a {\em loop} if $r^n(l)=d^n(l)$, 
and the vertex $r^n(l)=d^n(l)$ is called the {\em base point} of the loop $l$. 
A loop $l=(e_1,\ldots,e_n)$ is 
said to be {\em without entrances} 
if $r^{-1}(r(e_k))=\{e_{k}\}$ for $k=1,\ldots,n$.

\begin{definition}[{\cite[Definition 5.4]{Ka4}}]
A topological graph $E$ is said to be {\em topologically free} 
if the set of base points of loops without entrances
has an empty interior.  
\end{definition}

\begin{definition}[{\cite[Definition 8.4]{Ka6}}]
A topological graph $E$ is said to be {\em generated by a loop} $l$ 
if $E^0$ is discrete and 
every negative orbits are in the form 
$(\mu,l,l,\ldots)\in E^\infty$ with some $\mu\in E^*$. 
\end{definition}

\begin{proposition}[{\cite[Theorem 8.12]{Ka6}}]\label{PropSimple}
For a topological graph $E$, 
the following conditions are equivalent: 
\begin{enumerate}
\rom
\item The $C^*$-algebra $\cO(E)$ is simple. 
\item $E$ is minimal and topologically free. 
\item $E$ is minimal and not generated by a loop. 
\end{enumerate}
\end{proposition}

We need the following lemma in the next section 
which is a consequence of the technical result \cite[Proposition 5.10]{Ka4}. 

\begin{lemma}\label{LemPI2}
Let $E=(E^0,E^1,d,r)$ be a topologically free topological graph, 
and $v_0$ be an element of $E^0$ with $\overline{\Orb^+(v_0)}=E^0$. 
For a non-zero positive element $x\in \cO(E)$, 
there exist $a\in \cO(E)$ 
and $f\in C_0(E^0)$ 
which is $1$ on some neighborhood $V_0$ of $v_0$ 
such that $\|a^*xa-t^0(f)\|<1/2$.  
\end{lemma}

\begin{proof}
Take a non-zero positive element $x\in \cO(E)$. 
Set $\e=\|\varPsi(x)\|/5>0$ where 
$\varPsi\colon \cO(E)\to \cO(E)^{\gamma}$ 
is the faithful conditional expectation defined by
$$\varPsi(x)=\int_{\T}\gamma_z(x)dz,$$
using the gauge action $\gamma$ on $\cO(E)$. 
Choose a non-zero positive element $y\in \cO(E)$ 
with $\|x-y\|<\e$ that is a finite sum of elements 
in the form $t^n(\xi)t^m(\eta)^*$. 
Since $E$ is topologically free, 
there exist $a_0\in \cO(E)$ 
and a positive function $f_0\in C_0(E^0)$ 
such that $\|a_0\|\leq 1$, 
$\|f_0\|=\|\varPsi(y)\|$ 
and $\|a_0^*ya_0-t^0(f_0)\|<\e$ 
by \cite[Proposition 5.10]{Ka4} (see \cite[Remark 5.11]{Ka4}). 
Since 
\[
\|f_0\|=\|\varPsi(y)\|>\|\varPsi(x)\|-\e=4\e, 
\]
there exists a non-empty open set $V$ 
such that $f_0(v)>4\e$ for $v\in V$. 
Since $\Orb^+(v_0)$ is dense in $E^0$, 
we can find $\mu\in E^n$ with $d^n(\mu)=v_0$ 
and $r^n(\mu)\in V$. 
Choose a neighborhood $U$ of $\mu\in E^n$ 
such that the restriction of $d^n$ to $U$ is injective, 
and $f_0(v)>4\e$ for all $v\in r^n(U)$. 
We can find $\xi\in C_c(U)\subset C_{d^n}(E^n)$ 
such that $\|\xi\|^2<(4\e)^{-1}$ and 
$f=\ip{\xi}{\pi_{r^n}(f_0)\xi}\in C_0(E^0)$ is $1$ 
on some neighborhood $V_0$ of $v_0$. 
We set $a=a_0t^n(\xi)\in \cO(E)$. 
Then we have 
\begin{align*}
\|a^*xa-t^0(f)\|
&=\|t^n(\xi)^*a_0^*xa_0t^n(\xi)-t^n(\xi)^*t^0(f_0)t^n(\xi)\|\\
&\leq\big\|a_0^*xa_0-t^0(f_0)\big\|\|t^n(\xi)\|^2\\
&\leq\big(\|a_0^*(x-y)a_0\|+\|a_0^*ya_0-t^0(f_0)\|\big)\|\xi\|^2\\
&<(\e+\e)(4\e)^{-1} \\
&=1/2
\end{align*}
We are done. 
\end{proof}

\section{Contracting topological graphs}\label{SecCont}

In this section, 
we define contracting topological graphs, 
and prove Theorem A. 

\begin{definition}
Let $n,m$ be positive integers, and set $k=\min\{n,m\}$. 
For two subsets $U\subset E^n$ and $U'\subset E^m$, 
we define $U\pitchfork U'\subset E^k$ by 
$U\pitchfork U'=(U|_k)\cap (U'|_k)$ 
where 
\[
U|_k=\{(e_1,e_2,\ldots,e_k)\in E^k\mid 
(e_1,e_2,\ldots,e_n)\in U\}
\]
and $U'|_k$ is defined similarly. 
\end{definition}

Note that for $U,U'\subset E^n$, 
$U\pitchfork U'=U\cap U'$. 
The following follows from \cite[Lemma 2.4]{Ka4} easily. 

\begin{lemma}\label{Lem:apart}
Let $n,m$ be positive integers, 
and $U\subset E^n$ and $U'\subset E^m$ 
be open sets satisfying $U\pitchfork U'=\emptyset$. 
Then for any $\xi\in C_c(U)\subset C_{d^n}(E^n)$ 
and $\eta\in C_c(U')\subset C_{d^m}(E^m)$, 
we have $t^n(\xi)^*t^m(\eta)=0$. 
\end{lemma}

\begin{definition}\label{DefCont}
Let $E=(E^0,E^1,d,r)$ be a topological graph. 
We say that a non-empty open subset $V$ of $E^0$ 
is a {\em contracting} open set 
if its closure $\overline{V}$ is compact 
and there exist non-empty open subsets 
$U_k\subset E^{n_k}$ for $k=1,2,\ldots,m$ 
with $n_k\geq 1$ satisfying 
\begin{enumerate}
\rom
\item $r^{n_k}(U_k)\subset V$ for $k=1,2,\ldots,m$, 
\item $U_k\pitchfork U_l=\emptyset$ for $k\neq l$, 
\item $\overline{V}\subsetneq\bigcup_{k=1}^m d^{n_k}(U_k)$. 
\end{enumerate}
\end{definition}

If an open subset $V$ of $E^0$ 
with compact closure 
satisfies $\overline{V}\subsetneq d^n((r^n)^{-1}(V))$ 
for some $n\geq 1$, 
then it is contracting. 
The converse is not true in general (see Example \ref{Ex1}). 

\begin{lemma}\label{LemPI1}
If a topological graph $E$ has a contracting open set $V$, 
then there exist non-zero elements $x,y\in \cO(E)$ 
satisfying $x^*xx=x$, $x^*xy=y$, $x^*y=0$ 
and $t^0(f)x=x$ for all $f\in C_0(E^0)$ 
which is $1$ on $V$. 
\end{lemma}

\begin{proof}
Take non-empty open subsets $U_k\subset E^{n_k}$ for $k=1,2,\ldots,m$ 
satisfying (i), (ii) and (iii) in Definition \ref{DefCont}. 
Since $\overline{V}\subsetneq\bigcup_{k=1}^m d^{n_k}(U_k)$, 
there exists $k_0\in\{1,2,\ldots,m\}$ 
with $d^{n_{k_0}}(U_{k_0})\setminus \overline{V}\neq \emptyset$. 
Then we can find non-empty open sets $U_{k_0}',U_0\subset U_{k_0}$ 
with $U_{k_0}'\cap U_0=\emptyset$ and 
$\overline{V}\subset\bigcup_{k\neq k_0}d^{n_k}(U_k)
\cup d^{n_{k_0}}(U'_{k_0})$. 
Replacing $U_{k_0}$ by $U'_{k_0}$ and setting $n_0=n_{k_0}$, 
we get non-empty open subsets $U_k\subset E^{n_k}$ 
for $k=0,1,2,\ldots,m$ satisfying 
\begin{enumerate}
\renewcommand{\labelenumi}{{\rm (\roman{enumi})'}}%
\renewcommand{\itemsep}{0pt}
\item $r^{n_k}(U_k)\subset V$ for $k=0,1,2,\ldots,m$, 
\item $U_k\pitchfork U_l=\emptyset$ 
for $k,l\in\{0,1,2,\ldots,m\}$ with $k\neq l$, 
\item $\overline{V}\subset\bigcup_{k=1}^m d^{n_k}(U_k)$. 
\end{enumerate}
By (iii)', 
we can find $\xi_k\in C_c(U_k)\subset C_{d^{n_k}}(E^{n_k})$ 
for $k=1,\ldots,m$ 
such that $g=\sum_{k=1}^m \ip{\xi_k}{\xi_k}\in C_0(E^0)$ 
is $1$ on $\overline{V}$. 
Set $x=\sum_{k=1}^mt^{n_k}(\xi_k)$. 
Then for $f\in C_0(E^0)$ 
which is $1$ on $V$, 
we have $t^0(f)x=x$ by (i)'. 
By (ii)' and Lemma \ref{Lem:apart}, 
we get 
\[
x^*x=\sum_{k,l=1}^mt^{n_k}(\xi_k)^*t^{n_l}(\xi_l)
=\sum_{k=1}^mt^{n_k}(\xi_k)^*t^{n_k}(\xi_k)
=t^0(g).
\]
Since $g$ is $1$ on $\overline{V}$, 
we have $x^*xx=x$. 
Take $\xi_0\in C_c(U_0)\subset C_{d^{n_0}}(E^{n_0})$ with $\xi_0\neq 0$, 
and set $y=t^{n_0}(\xi_0)$. 
Then it is easy to check $x^*xy=y$ and $x^*y=0$. 
We are done. 
\end{proof}

\begin{lemma}\label{LemPI0}
If a topological graph $E$ is minimal 
and has a contracting open set, 
then the \Ca $\cO(E)$ is simple 
and has an infinite projection. 
\end{lemma}

\begin{proof}
Let $E$ be a minimal topological graph 
having a contracting open set. 
It is easy to see 
that a topological graph generated by a loop 
does not have a contracting open set. 
Hence by Proposition \ref{PropSimple}, 
the \Ca $\cO(E)$ is simple. 
By Lemma \ref{LemPI1}, 
there exists $x\in\cO(E)$ with $x^*xx=x$ 
and $x^*x\neq xx^*$. 
Such an element is called a scaling element in \cite{BC}, 
and a simple \Ca containing a scaling element 
has an infinite projection 
(see \cite[Proposition 4.2]{Ka2}). 
Hence the \Ca $\cO(E)$ has an infinite projection. 
\end{proof}

\begin{remark}
There may be a good chance to show that $\cO(E)$ is purely infinite 
under the assumption of Lemma \ref{LemPI0}, 
though the author could not prove it. 
\end{remark}

\begin{definition}\label{Def(PI)}
We say that a topological graph $E$ is 
{\em contracting at} $v_0\in E^0$ 
if $\overline{\Orb^+(v_0)}=E^0$, 
and any neighborhood $V_0$ of $v_0$ 
contains a contracting open set $V\subset V_0$. 
We simply say that $E$ is {\em contracting} 
if $E$ is contracting at some $v_0\in E^0$. 
\end{definition}

Now we show Theorem A which says that 
for a minimal and contracting topological graph $E$, 
the \Ca $\cO(E)$ is simple and purely infinite. 

\begin{proof}[Proof of Theorem A]
Let $E$ be a minimal topological graph 
which is contracting at $v_0$. 
By Lemma \ref{LemPI0}, 
the \Ca $\cO(E)$ is simple 
and has an infinite projection $p$. 
To show that $\cO(E)$ is purely infinite, 
it suffices to see that 
for each non-zero positive element $x_0\in \cO(E)$, 
the hereditary subalgebra $\overline{x_0\cO(E)x_0}$ of $\cO(E)$ 
generated by $x_0$ has a projection which is equivalent to 
the infinite projection $p$. 

Take a non-zero positive element $x_0\in \cO(E)$. 
By Lemma \ref{LemPI0} and Proposition \ref{PropSimple}, 
$E$ is topologically free. 
Hence by Lemma \ref{LemPI2}, 
there exists $a\in \cO(E)$ 
and $f\in C_0(E^0)$ 
which is $1$ on some neighborhood $V_0$ of $v_0$ 
such that $\|a^*x_0a-t^0(f)\|<1/2$.  
Since there exists a contracting open set $V\subset V_0$ 
and $f$ is $1$ on $V$, 
there exist non-zero elements $x,y\in \cO(E)$ 
satisfying $x^*xx=x$, $x^*xy=y$, $x^*y=0$ 
and $t^0(f)x=x$ by Lemma \ref{LemPI1}. 
Since $y^*y\neq 0$ and $\cO(E)$ is simple, 
we can find $b_1,\ldots,b_l\in \cO(E)$ such that 
\[
\sum_{k=1}^lb_k^*y^*yb_k=p.
\]
Set $b=\sum_{k=1}^lx^kyb_k\in \cO(E)$. 
Then we have $b^*t^0(f)b=b^*b=p$. 
This implies $\|b\|=1$. 
Therefore we get 
\[
\|b^*a^*x_0ab-p\|=
\|b^*(a^*x_0a-t^0(f))b\|<1/2. 
\]
Let $\chi$ be the characteristic function of 
an open interval $(1/2,3/2)$. 
Then $\chi(b^*a^*x_0ab)$ is a projection which is equivalent to $p$. 
The projection $q=\chi(x_0^{1/2}abb^*ax_0^{1/2})\in \overline{x_0\cO(E)x_0}$ 
is equivalent to $\chi(b^*a^*x_0ab)$, hence to 
the infinite projection $p$. 
We are done. 
\end{proof}

\begin{remark}
The author does not know 
whether the converse of Theorem A is true or not. 
\end{remark}

\section{Examples and remarks}\label{SecRem}

In this section, 
we give examples of minimal contracting topological graphs, 
and state some problems and remarks. 

A minimal topologically free discrete graph $E$ is contracting 
if and only if it has a loop. 
Thus for a discrete graph $E$, 
the \Ca $\cO(E)$ is purely infinite 
if $E$ is minimal, topologically free, and has a loop. 
The author does not know 
whether this is true for a topological graph $E$. 

\begin{problem}\label{Prob4} 
Suppose that a topological graph $E$ is minimal, topologically free. 
Is $\cO(E)$ purely infinite when $E$ has a loop? 
\end{problem}

It is known that for a minimal topologically free discrete graph $E$, 
having a loop is not only sufficient but also necessary 
for the \Ca $\cO(E)$ to be purely infinite 
(see the remark after Proposition \ref{ET:minimal}). 
However the following example shows that 
this is not the case for general topological graphs $E$. 

\begin{example}\label{Ex1}
Set $E^0=\R$ and $E^1=\R\times\{0,1\}$. 
We define two maps $d,r\colon E^1\to E^0$ 
by $d(x,0)=d(x,1)=x$, 
$r(x,0)=x+\sqrt{2}$ and $r(x,1)=x-1$. 
Then the topological graph $E=(E^0,E^1,d,r)$ has no loop, 
but $\cO(E)$ is simple and purely infinite. 
In fact, all $v_0\in E^0$ satisfies $\overline{\Orb^+(v_0)}=E^0$ 
and $E$ is contracting at all $v_0\in E^0$ 
because all non-empty open subsets of $E^0$ with compact closure 
are contracting. 
Note that if $V\subset E^0$ is small enough, 
there is no $n\geq 1$ 
with $\overline{V}\subset d^n((r^n)^{-1}(V))$. 
\end{example}

\begin{remark}
The \Ca $\cO(E)$ in the above example 
is naturally isomorphic to 
the crossed product $\cO_2\rtimes_{\alpha}\R$, 
where the action 
$\alpha\colon\R\curvearrowright \cO_2$ is 
defined by 
$\alpha_t(S_0)=e^{\sqrt{2}ti}S_0$, 
$\alpha_t(S_1)=e^{-ti}S_1$ 
for $t\in\R$ and 
for the standard generators $S_0,S_1$ 
of the Cuntz algebra $\cO_2$ (see \cite{KaK}). 
By computing the $K$-groups, 
we see that $\cO(E)\cong\cO_2\otimes\cK$ 
where $\cK$ is the \Ca of the all compact operators on 
the separable infinite dimensional Hilbert space. 
A similar example can be obtained by replacing $\R$ 
to $\R/2\pi\Z$. 
In that case, the \Ca is isomorphic to 
$\cO_2\rtimes_{\alpha}\Z\cong\cO_2$. 
\end{remark}

Although the topological graph in Example \ref{Ex1} has no loop, 
it has many ``approximate loops''. 

\begin{definition}
An {\em approximate loop} is a sequence $\{e_n\}_{n\in\N}$ in $E^1$ 
such that $d(e_n)=r(e_{n-1})$ for $n=1,2,\ldots,$ 
and that for every neighborhood $V$ of $d(e_0)$ 
infinitely many $n$ satisfy $r(e_n)\in V$. 
The vertex $v=d(e_0)\in E^0$ is called the {\em base point} 
of the approximate loop $e$. 
\end{definition}

For a loop $l=(e_1,e_2,\ldots,e_k)$ based on $v_0\in E^0$, 
the sequence $\{e'_n\}_{n\in\N}$ defined by $e'_{kl-j}=e_j$ 
is an approximate loop based on $v_0$. 
When $E^0$ is discrete, 
every approximate loops have this form. 

\begin{problem}
Suppose that a topological graph $E$ is minimal and topologically free. 
Are the following two conditions for $v_0\in E^0$ equivalent? 
\begin{enumerate}
\rom
\item $v_0$ is the base points of two distinct approximate loops, 
\item $E$ is contracting at $v_0\in E^0$. 
\end{enumerate}
\end{problem}

The following fact may be worth to remark. 

\begin{lemma}
Let $E$ be a minimal topological graph. 
If $v_0\in E^0$ is a base point of an approximate loop, 
then we have $\overline{\Orb^+(v_0)}=E^0$. 
\end{lemma}

\begin{proof}
Since $v_0$ is a base point of an approximate loop, 
either $v_0$ is a base point of a loop 
or $v_0$ is not isolated in $\Orb^+(v_0)$. 
In both cases, 
\cite[Proposition 4.4]{Ka6} implies 
that the closed set $\overline{\Orb^+(v_0)}$ is invariant 
and hence there exists a negative orbit $\mu$ of $v_0$ with 
$\Orb(v_0,\mu)\subset \overline{\Orb^+(v_0)}$ 
(see \cite[Proposition 4.11]{Ka6}). 
Since $E$ is minimal, 
we have $\overline{\Orb^+(v_0)}=E^0$. 
\end{proof}

For the topological graph $E$ in Example \ref{Ex1}, 
every vertices are base points of 
infinitely many approximate loops. 

\begin{problem}
Suppose that a topological graph $E$ is minimal, topologically free. 
If there exists $v_0\in E^0$ 
which is the base points of two distinct infinite loops, 
then is $\cO(E)$ purely infinite? 
How about the converse? 
\end{problem}

The existence of one approximate loop is not sufficient 
for pure infiniteness of $\cO(E)$ as the following example shows. 

\begin{example}
Define a topological graph $E=(E^0,E^1,d,r)$ 
by $E^0=\T$, $E^1=\T$, 
$d(x)=x$ and $r(x)=e^{2\pi i\theta}x$ 
where $\theta\in\R\setminus\Q$. 
Then $E$ is minimal, topologically free, 
and every vertices of $E$ are base points 
of approximate loops. 
However $\cO(E)$ is not purely infinite. 
In fact, $\cO(E)$ is isomorphic to the irrational rotation 
algebra $A_\theta$. 
Note that for each $v_0\in E^0$, 
there exists only one approximate loop whose base point is $v_0$. 
\end{example}

We finish this section by giving 
a method to construct a topological graph $E$ 
so that the \Ca $\cO(E)$ is a Kirchberg algebra 
with a given pair of countable abelian groups as its $K$-groups. 
This construction is a special case of the construction in \cite{Ku}. 

Let $G_0$ and $G_1$ be countable abelian groups. 
There exists a locally compact second countable space $X$ 
such that $K_i(C_0(X))\cong G_i$ for $i=0,1$ 
(for a proof, see \cite[Corollary 23.10.3]{B}). 
Let us choose a sequence $\{v_n\}_{n\in\N}$ in $X$
so that for all non-empty open set $V\subset X$, 
the set $\{n\in\N\mid v_n\in V\}$ is infinite. 
Let us set $E^0=X$, $E^1=X\times\N$ 
and define two maps $d,r$ from $E^1$ to $E^0$ 
by $d(x,n)=x$ and $r(x,n)=v_n$ 
for $(x,n)\in E^1$. 
Then the \Ca $\cO(E)$ arising from 
the topological graph $E=(E^0,E^1,d,r)$ 
satisfies the following. 

\begin{proposition}
The \Ca $\cO(E)$ is a Kirchberg algebra 
with $K_i(\cO(E))\cong G_i$ for $i=0,1$. 
\end{proposition}

\begin{proof}
It is clear to see that 
$E$ is minimal, 
and contracting at $v_n\in E^0$ for every $n$. 
Hence $\cO(E)$ is a Kirchberg algebra by Corollary B. 
It is also clear to see $\s{E}{0}{rg}=\emptyset$. 
Now Proposition \ref{Kgroup} shows 
$K_i(\cO(E))\cong K_i(C_0(E^0))\cong G_i$ for $i=0,1$. 
\end{proof}

By this proposition, 
we can see that the \Ca $\cO(E)$ does not depend on 
the choices of the sequence $\{v_n\}_{n\in\N}$. 
When $X$ is one point or a closed interval $[0,1]$, 
the \Ca $\cO(E)$ constructed as above is isomorphic to 
the Cuntz algebra $\cO_\infty$. 
When $X$ is a half-open interval $(0,1]$, 
the \Ca $\cO(E)$ is isomorphic to $\cO_2\otimes\cK$. 

By this construction, 
we get all non-unital Kirchberg algebras, 
but do not get all unital ones. 
The \Ca $\cO(E)$ is unital 
if and only if $X$ is compact (see \cite[Proposition 7.1]{Ka5}), 
and in this case the isomorphism $K_0(C_0(X))\cong K_0(\cO(E))$ 
in the above proof sends $[1_{C_0(X)}]\in K_0(C_0(X))$ 
to $[1_{\cO(E)}]\in K_0(\cO(E))$. 
Since $[1_{C_0(X)}]\in K_0(C_0(X))$ satisfies the property 
that $\Z[1_{C_0(X)}]$ is a direct summand of $K_0(C_0(X))$, 
so does $[1_{\cO(E)}]\in K_0(\cO(E))$ for a topological graph $E$ 
constructed from the compact space $X$ as above. 
Therefore by this construction 
we cannot get unital Kirchberg algebras $A$ such that 
$\Z[1_{A}]$ is not a direct summand of $K_0(A)$, 
such as the Cuntz algebras $\cO_n$ for $n<\infty$ 
or unital Kirchberg algebras $A$ 
with $K_0(A)\cong \Q$ or $[1_A]=0$. 
For a method to get all unital Kirchberg algebras 
from topological graphs, 
see Section \ref{SecKirch}.

\section{Topological graphs $E_{n,m}$}\label{SecEnm}

In this and the next sections, 
we give a way to construct a minimal contracting 
topological graph such that 
one can control the $K$-groups 
of the associated \CA s. 

Take $n\in\Z_+=\{1,2,\ldots\}$ and $m\in\Z$. 
We define a topological graph 
$E_{n,m}=(E^0,E^1,d,r)$ 
by $E^0=E^1=\T$, $d(z)=z^n$ and $r(z)=z^m$. 
We first compute the $K$-groups of the \Ca $\cO(E_{n,m})$ 
using Proposition \ref{Kgroup}. 
To this end, 
we need to know the subset $\s{E}{0}{rg}\subset E^0$ 
and the map $[\pi_r]$. 

\begin{lemma}\label{Enmrg}
When $m=0$ we have $\s{E}{0}{rg}=\emptyset$, 
and when $m\neq 0$ we have $\s{E}{0}{rg}=E^0$. 
\end{lemma}

\begin{proof}
When $m=0$, the image of $r$ is $\{1\}$. 
Hence $\s{E}{0}{rg}=\emptyset$. 
When $m\neq 0$, the map $r$ is a surjective proper map. 
Hence $\s{E}{0}{rg}=E^0$. 
\end{proof}

When $m\neq 0$, 
both the $K_0$-group and the $K_1$-group 
of $C_0(\s{E}{0}{rg})=C_0(E^0)=C(\T)$ are isomorphic to $\Z$ 
whose generators are given by $[1]$ and $[u]$ respectively, 
using the unit $1$ and the generating unitary $u$ 
of $C(\T)$. 
Thus when $m\neq 0$ 
the homomorphisms $[\pi_r]\colon K_i(C_0(\s{E}{0}{rg}))\to K_i(C_0(E^0))$ 
for $i=0,1$ 
are determined by the images of $[1]$ and $[u]$. 

For each $k\in\Z$, 
we define $\xi_k\in C_d(E^1)=C(\T)$ 
by $\xi_k(z)=z^k/\sqrt{n}$. 
Then it is not difficult to see $C_d(E^1)=\cspa\{\xi_k\mid k\in\Z\}$ 
and 
\[
\ip{\xi_{k'}}{\xi_k}=
\begin{cases}
u^l & \mbox{if $k-k'=nl$ for $l\in\Z$,}\\
0 & \mbox{if $k-k'\notin n\Z$.}
\end{cases}
\]
One can check the following by direct computations. 

\begin{lemma}\label{Enm1}
For $k,k'\in\Z$, 
we define $j,j'\in\{0,1,\ldots,n-1\}$ 
and $l,l'\in\Z$ 
by $k=j+nl$ and $k'=j'+nl'$. 
Then the map 
\[
\cK(C_d(E^1))\ni\theta_{\xi_k,\xi_{k'}}\mapsto 
e_{j,j'}\otimes u^{l-l'}\in \M_n\otimes C(\T)
\]
defines an isomorphism, 
where $\{e_{j,j'}\}_{j,j'=0}^{n-1}$ are 
the matrix units of the $n\times n$ matrix algebra $\M_n$. 
\end{lemma}

By Lemma \ref{Enm1}, 
the \Ca $\cK(C_d(E^1))$ has the unit 
which we denote by $I\in \cK(C_d(E^1))$. 
Define a unitary $U\in \cK(C_d(E^1))$ by 
\[
U=\sum_{j=0}^{n-1}\theta_{\xi_{j+1},\xi_j}\in \cK(C_d(E^1)).
\]
Then we have $U\xi_k=\xi_{k+1}$ for all $k\in\Z$. 
The image of $U$ under the isomorphism in Lemma \ref{Enm1} is 
\[
\sum_{j=0}^{n-2}e_{j+1,j}\otimes 1+e_{0,n-1}\otimes u
\in \M_n\otimes C(\T),
\] 
which is a generator of $K_1(\M_n\otimes C(\T))\cong\Z$. 
One can see that $[I]\in K_0\big(\cK(C_d(E^1))\big)$ and 
$[U]\in K_1(\cK(C_d(E^1)))$ correspond to 
$n[1]\in K_0(C_0(E^0))$ and $[u]\in K_1(C_0(E^0))$, respectively, 
under the isomorphism $K_i\big(\cK(C_d(E^1))\big)\cong K_i(C_0(E^0))$ 
determined by the full Hilbert $C_0(E^0)$-module $C_d(E^1)$. 

\begin{lemma}\label{Enm3}
When $m\neq 0$, 
the map $[\pi_r]\colon K_0(C(\s{E}{0}{rg}))\to K_0(C(E^0))$ 
sends $[1]$ to $n[1]$, 
and the map $[\pi_r]\colon K_1(C(\s{E}{0}{rg}))\to K_1(C(E^0))$ 
sends $[u]$ to $m[u]$. 
\end{lemma}

\begin{proof}
This follows from the facts $\pi_r(1)=I$ and $\pi_r(u)=U^m$, 
and the above computations. 
\end{proof}

Now we can compute the $K$-groups of the \Ca $\cO(E_{n,m})$. 

\begin{proposition}\label{KOEnm}
The pair of the $K$-groups $\big(K_0(\cO(E_{n,m})),K_1(\cO(E_{n,m}))\big)$ 
is as follows: 
\begin{center}
\begin{tabular}{c|c|c|c|}
& $m=0$ & $m=1$ & $m\neq 0,1$ \\
\hline
\multicolumn{1}{@{}c@{}|}{%
$\begin{array}{c}
\raisebox{3pt}[16pt]{$n=1$}\\ \hline 
\raisebox{3pt}[16pt]{$n\geq 2$}
\end{array}$}
& $\ (\Z,\ \Z)\ $ & 
\multicolumn{1}{@{}c@{}|}{%
$\begin{array}{c}
\raisebox{3pt}[16pt]{$(\Z\oplus\Z,\ \Z\oplus\Z)$}\\ \hline 
\raisebox{3pt}[16pt]{$\ (\Z\oplus\Z/(n-1)\Z,\ \Z)\ $}
\end{array}$}
& 
\multicolumn{1}{@{}c@{}|}{%
$\begin{array}{c}
\raisebox{3pt}[16pt]{$(\Z,\ \Z\oplus\Z/|m-1|\Z)$}\\ \hline 
\raisebox{3pt}[16pt]{$\ (\Z/(n-1)\Z,\ \Z/|m-1|\Z)\ $}
\end{array}$}\\
\hline 
\end{tabular}
\end{center}
The element $[1]\in K_0(\cO(E_{n,m}))$ corresponds 
to $1$ when $m\neq 1$, 
and to $(0,1)$ when $m=1$. 
\end{proposition}

\begin{proof}
The case $m=0$ follows from 
Proposition \ref{Kgroup} and Lemma \ref{Enmrg}, 
and the case $m\neq 0$ follows from 
Proposition \ref{Kgroup} and Lemma \ref{Enm3}. 
\end{proof}

Note that the \Ca $\cO(E_{1,1})$ is isomorphic to $C(\T\times\T)$. 
Next we determine 
for which $n,m$ the topological graph $E_{n,m}$ is minimal and contracting. 

\begin{proposition}\label{EnmMin}
When $m\in n\Z$, the topological graph $E_{n,m}$ is not minimal. 

When $m\notin n\Z$, 
$E_{n,m}$ is minimal and contracting, 
hence $\cO(E_{n,m})$ is a Kirchberg algebra. 
\end{proposition}

\begin{proof}
The infinite path $1^\infty=(1,1,\ldots)\in E^\infty$ 
is a negative orbit of $1\in E^0$, 
and we have $\Orb(1,1^\infty)=\{1\}$ when $m\in n\Z$. 
Hence in this case, $E_{n,m}$ is not minimal. 

Suppose $m\notin n\Z$. 
Let $p,q\in\Z$ be the relatively prime integers 
with $p>0$ and $q/p=m/n$. 
From $m\notin n\Z$ we get $p\geq 2$ and $q\neq 0$. 
For $e^{2\pi i \theta}\in E^0$, 
we have 
\begin{align*}
r\big(d^{-1}(e^{2\pi i \theta})\big)
&=r\big(\{e^{2\pi i (\theta+j)/n}\mid j=0,1,\ldots,n-1\}\big)\\
&=\{e^{2\pi i (\theta+j)m/n}\mid j=0,1,\ldots,n-1\}\\
&=\{e^{2\pi i (\theta+j)q/p}\mid j=0,1,\ldots,p-1\}\\
&=\{e^{2\pi i (\theta q+j)/p}\mid j=0,1,\ldots,p-1\},
\end{align*}
where in the last equality we use the fact that 
$q$ is relatively prime to $p$. 
By induction on $k\in\N$, 
we can show 
\[
r^k\big((d^k)^{-1}(e^{2\pi i \theta})\big)
=\big\{e^{2\pi i (\theta q^k+j)/p^k}\ \big|\ j=0,1,\ldots,p^k-1\big\}.
\]
Since $p\geq 2$, 
we have $\overline{\Orb^+(e^{2\pi i \theta})}=E^0$ 
for every $e^{2\pi i \theta}\in E^0$. 
Hence $E_{n,m}$ is minimal. 

We will show that $E_{n,m}$ is contracting at $1\in E^0$. 
To do so, it suffices to show that for any $\e>0$ 
there exists $k$ such that $d^k\big((r^k)^{-1}(V_\e)\big)=E^0$ 
where $V_\e=\{e^{2\pi i \theta}\mid -\e<\theta<\e\}$ 
is a neighborhood of $1$. 
Similarly as above, we see 
\[
d^k\big((r^k)^{-1}(e^{2\pi i \theta})\big)
=\big\{e^{2\pi i (\theta p^k+j)/q^k}\ \big|\ j=0,1,\ldots,|q|^k-1\big\}. 
\]
Hence 
\[
d^k\big((r^k)^{-1}(V_\e)\big)
=\bigcup_{j=0}^{|q|^k-1} \{e^{2\pi i (\theta p^k+j)/q^k}\mid -\e<\theta<\e\}. 
\]
Thus $d^k\big((r^k)^{-1}(V_\e)\big)=E^0$ 
when $k$ is large enough that $2\e p^k>1$. 
We are done. 
\end{proof}

\section{Topological graphs $E\times_{n,m}\T$}\label{SecEnmT}

Let $E=(E^0,E^1,d,r)$ be a discrete graph. 
Let $n\colon E^1\to\Z_+$ and $m\colon E^1\to\Z$ be two maps. 
We define two continuous maps 
$\td,\tr\colon E^1\times\T\to E^0\times\T$ 
by $\td(e,z)=(d(e),z^{n(e)})$ and 
$\tr(e,z)=(r(e),z^{m(e)})$ for $(e,z)\in E^1\times\T$. 
Since $n(e)\geq 1$ for all $e\in E^1$, 
the map $\td$ is locally homeomorphic. 
Hence we get 
a topological graph $(E^0\times\T, E^1\times\T,\td,\tr)$. 

\begin{definition}
We denote by $E\times_{n,m}\T$ 
the topological graph $(E^0\times\T, E^1\times\T,\td,\tr)$ 
defined as above. 
\end{definition}

\begin{example}
Let $E=(\{v\},\{e\},d,r)$ be 
the graph consisting of one vertex and one loop. 
For $n\in\Z_+$ and $m\in\Z$, we denote two maps 
$e\mapsto n$ and $e\mapsto m$ by the same symbols $n,m$. 
Then the topological graph $E\times_{n,m}\T$ 
is isomorphic to the topological graph $E_{n,m}$ constructed 
in the previous section. 
\end{example}

\begin{example}
Let $E=(\{v,w\},\{e\},d,r)$ be 
the graph consisting of two vertices and one edge 
with $d(e)=v$ and $r(e)=w$. 
For $n\in\Z_+$ and $m\in\Z$, we denote two maps 
$e\mapsto n$ and $e\mapsto m$ by the same symbols $n,m$. 
Then one can verify that 
the \Ca $\cO(E\times_{n,m}\T)$ is isomorphic to 
$\M_{n+1}\otimes C(\T)$ when $m\neq 0$, 
and to $\big(\M_{n+1}\otimes C(\T)\big)\oplus C(\T)$ 
when $m=0$ 
(see Remark \ref{lastRemark}). 
\end{example}

Take a discrete graph $E=(E^0,E^1,d,r)$ and 
two maps $n\colon E^1\to\Z_+$, $m\colon E^1\to\Z$, 
and fix them. 
We first compute the $K$-groups of the \Ca $\cO(E\times_{n,m}\T)$. 

\begin{definition}
We define $E^0_{m}\subset E^0$ by 
\[
E^0_{m}
=\{v\in E^0\mid 0<|r^{-1}(v)\setminus m^{-1}(0)|<\infty\}.
\]
\end{definition}

If $m^{-1}(0)=\emptyset$, 
then we have $E^0_{m}=E^0_{\rs{rg}}$. 
In general, $E^0_{m}$ need not contain nor be contained in 
$E^0_{\rs{rg}}$. 

\begin{lemma}\label{Lem:EnmTrg}
We have 
\[
(E\times_{n,m}\T)^0_{\rs{rg}}
=\big((E^0_{\rs{rg}}\cap E^0_{m})\times\T\big)\cup 
\big((E^0_{\rs{sg}}\cap E^0_{m})\times (\T\setminus\{1\})\big). 
\]
\end{lemma}

\begin{proof}
It is routine to check that 
for $v\in E^0_{\rs{rg}}\cap E^0_{m}$ 
the restriction of $\tr$ to 
$\tr^{-1}(\{v\}\times\T)=r^{-1}(v)\times\T$ 
is a proper surjection onto $\{v\}\times\T$, 
and that for $v\in E^0_{\rs{sg}}\cap E^0_{m}$ 
the restriction of $\tr$ to 
$\tr^{-1}(\{v\}\times(\T\setminus\{1\}))$ 
is a proper surjection onto $\{v\}\times(\T\setminus\{1\})$. 
This proves the inclusion 
\[
(E\times_{n,m}\T)^0_{\rs{rg}}
\supset \big((E^0_{\rs{rg}}\cap E^0_{m})\times\T\big)\cup 
\big((E^0_{\rs{sg}}\cap E^0_{m})\times (\T\setminus\{1\})\big). 
\]
For $v\in E^0_{\rs{sg}}\cap E^0_{m}$, 
we have $(v,1)\notin (E\times_{n,m}\T)^0_{\rs{rg}}$ 
because $\tr^{-1}(v,1)=r^{-1}(v)\times\{1\}$ is not compact. 
Finally take $v\in E^0\setminus E^0_{m}$. 
The set $r^{-1}(v)\setminus m^{-1}(0)$ is either empty or infinite. 
If $r^{-1}(v)\setminus m^{-1}(0)=\emptyset$ 
then the image of the restriction of $\tr$ 
to $\tr^{-1}(\{v\}\times\T)$ is $\{v\}\times\{1\}$, 
and if $|r^{-1}(v)\setminus m^{-1}(0)|=\infty$ 
then $\tr^{-1}(v,z)$ is not compact for all $z\in\T$. 
Hence $(v,z)\notin (E\times_{n,m}\T)^0_{\rs{rg}}$ 
for $(v,z)\in (E^0\setminus E^0_{m})\times\T$. 
This completes the proof. 
\end{proof}

For $v\in E^0$, 
let $p_v$ and $u_v$ 
be the unit and the generating unitary 
of $C(\{v\}\times\T)\subset C_0((E\times_{n,m}\T)^0)$, 
respectively. 
Then 
\begin{align*}
\Z^{E^0}\ni (k_v)_{v\in E^0}&\mapsto 
\sum_{v\in E^0}k_v[p_v]\in K_0\big(C_0((E\times_{n,m}\T)^0)\big),\\
\Z^{E^0}\ni (k_v)_{v\in E^0}&\mapsto 
\sum_{v\in E^0}k_v[u_v]\in K_1\big(C_0((E\times_{n,m}\T)^0)\big)
\end{align*}
are isomorphisms. 
By similar isomorphisms, 
we have 
\[
\Z^{E^0_{\rs{rg}}\cap E^0_{m}}\cong 
K_0\big(C_0((E\times_{n,m}\T)^0_{\rs{rg}})\big),\quad
\Z^{E^0_{m}}\cong 
K_1\big(C_0((E\times_{n,m}\T)^0_{\rs{rg}})\big)
\]
by Lemma \ref{Lem:EnmTrg}.

\begin{lemma}\label{EnmTK}
For $w\in E^0_{\rs{rg}}\cap E^0_{m}$ 
we have $[\pi_{\tr}]([p_w])=\sum_{e\in r^{-1}(w)}n(e)[p_{d(e)}]$, 
and for $w\in E^0_{m}$ 
we have $[\pi_{\tr}]([u_w])=\sum_{e\in r^{-1}(w)}m(e)[u_{d(e)}]$. 
\end{lemma}

\begin{proof}
This follows from Lemma \ref{Enm3}. 
\end{proof}

Let $I_0\colon \Z^{E^0_{\rs{rg}}\cap E^0_{m}}\to \Z^{E^0}$ 
and $I_1\colon \Z^{E^0_{m}}\to \Z^{E^0}$ 
be the embedding maps. 
Define an $E^0\times (E^0_{\rs{rg}}\cap E^0_{m})$-matrix $N$ 
by $N_{v,w}=\sum_{e\in d^{-1}(v)\cap r^{-1}(w)}n(e)$ 
for $v\in E^0$ and $w\in E^0_{\rs{rg}}\cap E^0_{m}$. 
Note that $N_{v,w}=0$ if there exists no $e\in E^1$ 
with $d(e)=v$ and $r(e)=w$. 
Similarly, 
we define an $E^0\times E^0_{m}$-matrix $M$ 
by $M_{v,w}=\sum_{e\in d^{-1}(v)\cap r^{-1}(w)}m(e)$. 

\begin{proposition}\label{Prop:K}
We have 
\begin{align*}
K_0\big(\cO(E\times_{n,m}\T)\big)
&\cong \coker (I_0-N) \oplus \ker (I_1-M),\\
K_1\big(\cO(E\times_{n,m}\T)\big)
&\cong \ker (I_0-N)\oplus \coker (I_1-M).
\end{align*}
\end{proposition}

\begin{proof}
This follows from Proposition \ref{Kgroup} and Lemma \ref{EnmTK}. 
\end{proof}

\begin{corollary}\label{Cor:unitK}
The \Ca $\cO(E\times_{n,m}\T)$ is unital 
if and only if $E^0$ is finite. 
In this case, the following holds. 
\begin{enumerate}
\rom 
\item the element $[1_{\cO(E\times_{n,m}\T)}]$ in $K_0(\cO(E\times_{n,m}\T))$ 
corresponds to the image of $(1,\ldots,1)\in \Z^{E^0}$ 
in $\coker (I_0-N)$ via the isomorphism in Proposition \ref{Prop:K}, 
\item both $K_0(\cO(E\times_{n,m}\T))$ and $K_1(\cO(E\times_{n,m}\T))$ 
are finitely generated, 
\item the rank of $K_1(\cO(E\times_{n,m}\T))$ 
does not exceed the one of $K_0(\cO(E\times_{n,m}\T))$. 
\end{enumerate}
\end{corollary}

Next we examine conditions on $E,n,m$ 
for the topological graph $E\times_{n,m}\T$ 
to be minimal or contracting. 
For $\mu=(e_1,\ldots,e_k)\in E^k$ with $k\geq 1$, 
we define $\T_\mu^k$ by 
\[
\T_\mu^k=\{
z=(z_1,\ldots,z_k)\in \T^k \mid 
\text{$z_i^{n(e_i)}=z_{i+1}^{m(e_{i+1})}$ 
for $i=1,\ldots,k-1$}\}.
\]
Define two maps 
$d_\mu,r_\mu\colon\T_\mu^k\to \T$ 
by $d_\mu(z)=z_k^{n(e_k)}$ and $r_\mu(z)=z_1^{m(e_1)}$ 
for $z=(z_1,\ldots,z_k)\in \T_\mu^k$. 
Then it is easy to see 
\[
(E\times_{n,m}\T)^k
=\{(\mu,z)\mid \mu\in E^k, z\in \T^k_\mu\},
\]
and two maps 
$\td^k,\tr^k\colon (E\times_{n,m}\T)^k\to E^0\times \T$ 
are expressed as $\td^k(\mu,z)=(d^k(\mu),d_\mu(z))$ 
and $\tr^k(\mu,z)=(r^k(\mu),r_\mu(z))$. 

\begin{definition}
For $e\in E^1$, we define $p(e)\in\Z_+$ by 
\[
p(e)=
\begin{cases}
1 & \text{if $m(e)=0$,}\\
n(e)/\big(n(e),|m(e)|\big) & \text{if $m(e)\neq 0$.} 
\end{cases}
\]
Recursively, 
for $\mu=(e,\nu)\in E^{k+1}$ with $e\in E^1$ and $\nu\in E^k$, 
we define $p(\mu)\in\Z_+$ by 
\[
p(\mu)=
\begin{cases}
1 & \text{if $m(e)=0$,}\\
n(e)p(\nu)/\big(n(e)p(\nu),|m(e)|\big) & \text{if $m(e)\neq 0$.} 
\end{cases}
\]
\end{definition}

\begin{lemma}\label{Lem:p(mu,nu)}
Let $\mu\in E^k$ and $\nu\in E^l$ 
with $k,l\geq 1$ and $d^k(\mu)=r^l(\nu)$, 
and set $\mu'=(\mu,\nu)\in E^{k+l}$. 
Then we have $p(\mu)\mid p(\mu')$. 
\end{lemma}

\begin{proof}
When $k=1$, 
the conclusion is clear from the definition of $p(\cdot)$. 
Now we can prove the statement by the induction on $k$ 
using the fact that 
$p\mid p'$ implies $np/(np,m)\mid np'/(np',m)$ 
for $n,m,p,p'\in\Z_+$. 
\end{proof}

\begin{lemma}\label{Lem:rede}
Let us take $\mu\in E^k$ with $k\geq 1$. 
For $z_0\in \T$, we get 
\[
r_\mu(d_\mu^{-1}(z_0))
=\{z_ke^{2\pi ij/p(\mu)}\mid j=0,1,\ldots,p(\mu)-1\}
\] 
for some element $z_k\in \T$. 
\end{lemma}

\begin{proof}
The proof goes by induction on $k$. 
The case for $k=1$ can be shown easily 
(see the proof of Proposition \ref{EnmMin}). 
Suppose that we have proved the statement for $k$, 
and take $\mu=(e,\nu)\in E^{k+1}$ 
with $e\in E^1$ and $\nu\in E^k$. 
By the assumption of the induction, 
we obtain 
\begin{align*}
r_\mu(d_\mu^{-1}(z_0))
&=r_{e}\big(d_{e}^{-1}(r_{\nu}(d_{\nu}^{-1}(z_0)))\big)\\
&=r_{e}\big(d_{e}^{-1}(
\{z_ke^{2\pi ij/p(\nu)}\mid j=0,1,\ldots,p(\nu)-1\})\big)
\end{align*}
for some element $z_k\in \T$.
Take $z_{k+1}'$ with $(z_{k+1}')^{n(e)}=z_k$, 
and set $z_{k+1}=(z_{k+1}')^{m(e)}$. 
Then 
\begin{align*}
r_\mu(d_\mu^{-1}(z_0))
&=r_{e}\big(
\{z_{k+1}'e^{2\pi ij/n(e)p(\nu)}\mid j=0,1,\ldots,n(e)p(\nu)-1\}
\big)\\
&=
\{z_{k+1}e^{2\pi ijm(\mu)/n(e)p(\nu)}\mid j=0,1,\ldots,n(e)p(\nu)-1\}\\
&=
\{z_{k+1}e^{2\pi ij/p(\mu)}\mid j=0,1,\ldots,p(\mu)-1\}.
\end{align*}
This completes the proof. 
\end{proof}

\begin{corollary}\label{Cor:rede}
For $\mu\in E^k$,  
$r_\mu(d_\mu^{-1}(1))
=\{e^{2\pi ij/p(\mu)}\mid j=0,1,\ldots,p(\mu)-1\}$. 
\end{corollary}

\begin{proof}
Follows from Lemma \ref{Lem:rede}. 
\end{proof}

For $\mu=(e_1,e_2,\ldots,e_n,\ldots)\in E^\infty$, 
we set 
\[
\T_\mu^\infty=\{
z=(z_1,\ldots,z_k,\ldots)\in \T^\infty \mid 
\text{$z_i^{n(e_i)}=z_{i+1}^{m(e_{i+1})}$ 
for $i=1,2,\ldots$}\}.
\]
Then we have 
$(E\times_{n,m}\T)^\infty
=\{(\mu,z)\mid \mu\in E^\infty, z\in \T^\infty_\mu\}$. 
For every $\mu\in E^k$ with $k\in\N\cup\{\infty\}$, 
we denote by $1^k$ the element $(1,\ldots,1)\in \T^k_\mu$. 

\begin{proposition}\label{ET:minimal}
The topological graph $E\times_{n,m}\T$ 
is minimal 
if and only if the following two conditions are satisfied:
\begin{enumerate}
\rom
\item for every $v_0\in E^0$, 
every negative orbit $\mu_0$ of $v_0$ and every $v\in E^0$, 
\[
\sup\{p(\mu)\mid \mu\in E^*, d^*(\mu)\in\Orb^-(v_0,\mu_0), r^*(\mu)=v\}
=\infty,
\]
\item for every $v_0\in E^0_{\rs{rg}}\setminus E^0_{m}$ 
and every $v\in E^0$, 
\[
\sup\{p(\mu)\mid \mu\in E^*, d^*(\mu)=v_0, r^*(\mu)=v\}
=\infty.
\]
\end{enumerate}
\end{proposition}

\begin{proof}
Suppose that the topological graph $E\times_{n,m}\T$ 
is minimal. 
Take $v_0\in E^0$, 
a negative orbit $\mu_0\in E^k$ of $v_0$ 
with $k\in\N\cup\{\infty\}$. 
Then 
$(\mu_0,1^k)\in (E\times_{n,m}\T)^k$ is a negative orbit 
of $(v_0,1)\in E^0\times\T$. 
For $v\in E^0$, 
we set $L_v=\{\mu\in E^*\mid d^*(\mu)\in\Orb^-(v_0,\mu_0), r^*(\mu)=v\}$. 
Then, 
\begin{align*}
\{z\in\T\mid (v,z)\in\Orb((v_0,1),(\mu_0,1^k))\}
&=\bigcup_{\mu\in L_v}r_\mu(d_\mu^{-1}(1))\\
&=\bigcup_{\mu\in L_v}\{e^{2\pi ij/p(\mu)}\mid j=0,1,\ldots,p(\mu)-1\}
\end{align*}
by Corollary \ref{Cor:rede}. 
Since $E\times_{n,m}\T$ is minimal, 
$\Orb((v_0,1),(\mu_0,1^k))$ is dense in $E^0\times\T$. 
This shows $\{p(\mu)\mid \mu\in L_v\}=\infty$ 
for all $v\in E^0$.
Next we take $v_0\in E^0_{\rs{rg}}\setminus E^0_{m}$. 
Then $(v_0,1)\in (E\times_{n,m}\T)^0$ 
is a negative orbit of $(v_0,1)$. 
Now in a similar way as above, 
we get 
\[
\sup\{p(\mu)\mid \mu\in E^*, d^*(\mu)=v_0, r^*(\mu)=v\}
=\infty.
\]
for all $v\in E^0$. 
Thus (i) and (ii) are satisfied. 

Conversely, suppose that two conditions (i) and (ii) 
are satisfied. 
Take $(v_0,z_0)\in E^0\times \T$, 
and a negative orbit $(\mu_0,z)\in (E\times_{n,m}\T)^k$ 
for $k\in\N\cup\{\infty\}$. 
Then either $\mu_0\in E^k$ 
is a negative orbit of $v_0\in E^0$, 
or $k<\infty$ and $d^k(\mu_0)\in E^0_{\rs{rg}}\setminus E^0_{m}$ 
by Lemma \ref{Lem:EnmTrg}. 
For both cases, 
the conditions (i) and (ii) 
together with Lemma \ref{Lem:rede} shows that 
$\Orb((v_0,z_0),(\mu_0,z))$ is dense in $E^0\times\T$. 
Thus $E\times_{n,m}\T$ is minimal. 
\end{proof}

By Proposition \ref{ET:minimal}, 
if $E\times_{n,m}\T$ is minimal, 
then so is $E$. 
For a minimal discrete graph $E$ 
which is not generated by a loop, 
there exists a beautiful dichotomy 
(see \cite[Theorem 18]{S}, \cite[Remark 2.16]{DT}); 
\begin{align*}
\text{$\cO(E)$ is simple and purely infinite}&\iff 
\text{$E$ has a loop,}\\
\text{$\cO(E)$ is a simple AF-algebra}&\iff 
\text{$E$ has no loop.}
\end{align*}
For minimal topological graphs in the form 
$E\times_{n,m}\T$, 
we also have a similar dichotomy 
(see the remark after Proposition \ref{ET:dico}). 

A {\em circle algebra} is a \Ca which is isomorphic to 
$A\otimes C(\T)$ for some 
finite dimensional \Ca $A$. 
Let's say a \Ca is an {\em AT-algebra} 
if it is isomorphic to an inductive limit 
of quotients of circle algebras. 
For separable \CA s, 
this notion coincides with the one 
in literatures (see \cite[Proposition 3.2.3]{RS}). 
The proof of the following lemma is done in 
Appendix \ref{app1}. 

\begin{lemma}\label{Lem:appendix}
Let $F^0\subset E^0$ and $F^1\subset E^1$ be 
finite subsets satisfying $d(F^1),r(F^1)\subset F^0$. 
If the finite graph $F=(F^0,F^1,d|_{F^1},r|_{F^1})$ 
has no loops, 
then the \Csa of $\cO(E\times_{n,m}\T)$ 
generated by $t^0(C(F^0\times\T))$ and $t^1(C(F^1\times\T))$ 
is isomorphic to a quotient of a circle algebra. 
\end{lemma}

\begin{proposition}\label{ATalg}
If $E$ has no loop, 
$\cO(E\times_{n,m}\T)$ is an AT-algebra. 
\end{proposition}

\begin{proof}
This follows from Lemma \ref{Lem:appendix}. 
\end{proof}

\begin{proposition}\label{ET:dico}
When $E\times_{n,m}\T$ is minimal, 
the following conditions are 
equivalent;
\begin{enumerate}
\rom
\item $E\times_{n,m}\T$ is contracting, 
\item $E\times_{n,m}\T$ has a loop, 
\item $E$ has a loop, 
\item $\cO(E\times_{n,m}\T)$ is simple and purely infinite. 
\end{enumerate}
If moreover both $E^0$ and $E^1$ are countable, 
then the above equivalent conditions 
are also equivalent to 
\begin{enumerate}
\setcounter{enumi}{3}
\renewcommand{\labelenumi}{{\rm (\roman{enumi})'}}
\renewcommand{\itemsep}{0pt}
\item $\cO(E\times_{n,m}\T)$ is a Kirchberg algebra. 
\end{enumerate}
\end{proposition}

\begin{proof}
By Theorem A, 
we get (i)$\Rightarrow$(iv). 
As explained before Corollary B, 
we have (iv)$\Leftrightarrow$(iv)' 
when $E^0$ and $E^1$ are countable. 
By Proposition \ref{ATalg}, 
we get (iv)$\Rightarrow$(iii). 
It is easy to see (ii)$\Leftrightarrow$(iii). 
We will show (iii)$\Rightarrow$(i). 
Take a loop $e_0$ in $E$. 
Let $v_0\in E^0$ be the base point of the loop $e_0$. 
The infinite path $e_0^\infty=(e_0,e_0,\ldots)$ is a negative orbit 
of $v_0\in E^0$. 
By Proposition \ref{ET:minimal}, 
we have 
\[
\sup\{p(\mu)\mid \mu\in E^*, 
d^*(\mu)\in \Orb^-(v_0,e_0^\infty),r^*(\mu)=v\}
=\infty 
\]
for every $v\in E^0$. 
For $\mu\in E^*$ with $d^*(\mu)\in\Orb^-(v_0,e_0^\infty)$ 
and $r^*(\mu)=v$, 
there exists $\nu\in E^*$ 
with $r^*(\nu)=d^*(\mu)$ and $d^*(\nu)=v_0$. 
Then $\mu'=(\mu,\nu)\in E^*$ 
satisfies $d^*(\mu')=v_0$, $r^*(\mu')=v$ 
and $p(\mu)\leq p(\mu')$ 
by Lemma \ref{Lem:p(mu,nu)}. 
Hence we get 
\[
\sup\{p(\mu)\mid \mu\in E^*, d^*(\mu)=v_0, r^*(\mu)=v\}
=\infty
\]
for every $v\in E^0$. 
This shows that $\Orb^+((v_0,1))$ 
is dense in $E^0\times\T$. 
We will show that $E\times_{n,m}\T$ 
is contracting at $(v_0,1)$. 
To do so, 
it suffices to show that for all $\e>0$, 
there exist $k\in\Z_+$ 
and an open subset $U\subset (E\times_{n,m}\T)^k$ 
such that 
$\tr^k(U)\subset \{v_0\}\times V_\e$ 
and $\td^k(U)=\{v_0\}\times\T$ 
where $V_\e=\{e^{2\pi i \theta}\mid -\e<\theta<\e\}$. 
For $\e>0$, 
there exists $\mu\in E^k$ 
with $k\geq 1$ and $d^k(\mu)=r^k(\mu)=v_0$ 
such that $1/p(\mu)<2\e$. 
By Lemma \ref{Lem:rede}, 
$r_\mu(d_\mu^{-1}(z_0))\cap V_\e\neq\emptyset$ 
for all $z_0\in \T$. 
Therefore the open subset 
$U=\{\mu\}\times r_\mu^{-1}(V_\e)
\subset (E\times_{n,m}\T)^k$ 
satisfies 
$\tr^k(U)\subset \{v_0\}\times V_\e$ 
and $\td^k(U)=\{v_0\}\times\T$. 
Thus $E\times_{n,m}\T$ is contracting at $(v_0,1)$. 
We are done. 
\end{proof}

By Propositions \ref{ATalg} and \ref{ET:dico}, 
we get a dichotomy;
\begin{align*}
\text{$\cO(E\times_{n,m}\T)$ is simple and purely infinite}&\iff 
\text{$E$ has a loop,}\\
\text{$\cO(E\times_{n,m}\T)$ is a simple AT-algebra}&\iff 
\text{$E$ has no loop,}
\end{align*}
when $E\times_{n,m}\T$ is minimal.

\section{Construction of all Kirchberg algebras}\label{SecKirch} 

In this section, 
we construct all Kirchberg algebras 
from topological graphs. 
We first show that the class of \Cas $\cO(E\times_{n,m}\T)$ 
considered in the previous section 
contains all non-unital Kirchberg algebras. 

We denote by $\Z^\infty$ the group of 
sequences $x=(x_k)_{k=1}^\infty$ such that 
$x_k=0$ for all but finite $k$, 
and by $M_\infty(\Z)$ the set of 
integer-valued matrices $S=(S_{k,l})_{k,l=1}^\infty$ 
such that for each $l$, 
we have $S_{k,l}=0$ for all but finite $k$. 
An element $S$ in $M_\infty(\Z)$ is considered as 
a group homomorphism from $\Z^\infty$ to itself 
by $(Sx)_k=\sum_{l=1}^{\infty}S_{k,l}x_l$ for 
$x=(x_k)_{k=1}^\infty\in\Z^\infty$. 
We denote by $I\in M_\infty(\Z)$ the identity homomorphism. 
The set $M_\infty(\N)$ is the subset of $M_\infty(\Z)$ 
consisting of the matrices whose entries are in $\N$. 

Take the pair $N\in M_\infty(\N)$ 
and $M\in M_\infty(\Z)$ 
such that $N_{k,l}=0$ implies $M_{k,l}=0$. 
We define a discrete graph $E=(E^0,E^1,d,r)$ 
by $E^0=\Z_+$, 
$E^1=\{(k,l)\in \Z_+\times \Z_+\mid N_{k,l}\geq 1\}$, 
$d\colon E^1\ni (k,l)\mapsto k\in E^0$ 
and $r\colon E^1\ni (k,l)\mapsto l\in E^0$. 
Define $n,m\colon E^1\to \Z$ 
by $n(k,l)=N_{k,l}$ and $m(k,l)=M_{k,l}$ for $(k,l)\in E^1$. 

\begin{definition}\label{DefENM}
Let $N\in M_\infty(\N)$ 
and $M\in M_\infty(\Z)$ 
such that $N_{k,l}=0$ implies $M_{k,l}=0$. 
We denote by $E_{N,M}$ the topological graph $E\times_{n,m}\T$ 
defined as above. 
\end{definition}

\begin{lemma}\label{ENMK}
If $M$ has no columns that are identically zero, 
then we have $E^0=\s{E}{0}{rg}=E^0_m$ and 
\begin{align*}
K_0\big(\cO(E_{N,M})\big)
&\cong \coker (I-N) \oplus \ker (I-M),\\
K_1\big(\cO(E_{N,M})\big)
&\cong \ker (I-N)\oplus \coker (I-M).
\end{align*}
\end{lemma}

\begin{proof}
The former is easy to see, and the latter follows from 
Proposition \ref{Prop:K}. 
\end{proof}

\begin{remark}\label{RemOENM}
Similarly as above, 
from two finite matrices $N\in M_K(\N)$ and $M\in M_K(\Z)$ 
such that $N_{k,l}=0$ implies $M_{k,l}=0$, 
we can construct the topological graph $E_{N,M}$. 
When $K=1$, 
$E_{N,M}$'s are examples considered in Section \ref{SecEnm}. 
In this case, 
the construction of the \Ca $\cO(E_{N,M})$ 
can be compared with the one of the Cuntz-Krieger algebras 
in \cite{CK}, 
and Deaconu considered a similar construction in \cite{D} 
and proved the above lemma in his situation. 
\end{remark}

The following lemma is inspired by \cite[Lemma 1.1]{S2}. 

\begin{lemma}\label{MatrixSt}
For countable abelian groups $G_0$ and $G_1$, 
there exist $N\in M_\infty(\N)$ 
and $M\in M_\infty(\Z)$ satisfying the following: 
\begin{enumerate}
\rom
\item $N_{k,l}=0$ implies $M_{k,l}=0$, 
\item $N_{k,k}\geq 2$, $M_{k,k}=1$ for all $k\in\Z_+$, 
\item for all $(k,l)\in \Z_+\times \Z_+$, 
there exist $(k_i,l_i)$ with $N_{k_i,l_i}\geq 1$ 
for $i=0,1,\ldots,m$ 
such that $k_0=k$, $l_i=k_{i+1}$ 
for $i=0,\ldots,m-1$, and $l_m=l$, 
\item there exist exact sequences 
\[
\begin{CD}
0 @>>> \Z^\infty @>I-N>> \Z^\infty @>>> G_0 @>>> 0,\\
0 @>>> \Z^\infty @>I-M>> \Z^\infty @>>> G_1 @>>> 0.
\end{CD}
\]
\end{enumerate}
\end{lemma}

\begin{proof}
Since $G_0$ and $G_1$ are countable, 
there exist injective homomorphisms 
$T,S\colon \Z^\infty\to \Z^\infty$ 
with $\coker T\cong G_0$ and $\coker S\cong G_1$. 
Let us set $T^+,T^-,|S|\in M_\infty(\N)$ 
by $T^+_{k,l}=\max\{T_{k,l},0\}$, $T^-_{k,l}=\max\{-T_{k,l},0\}$ and 
$|S|_{k,l}=|S_{k,l}|$ 
for $k,l\in\Z_+$. 
We have $T=T^+-T^-$. 
Define $X\in M_\infty(\N)$ by 
\[
X_{k,l}=\begin{cases}
1 & \text{if }|k-l|\leq 1,\\
0 & \text{if }|k-l|\geq 2.
\end{cases}
\]
Define $N\in M_2(M_\infty(\N))$ and $M\in M_2(M_\infty(\Z))$ by
\[
N=\mat{2I}{T^++|S|+X}{I}{I+T^-+|S|+X},\quad
M=\mat{I}{S}{I}{I}.
\]
Let $I^{(2)}\in M_2(M_\infty(\Z))$ 
be the identity homomorphism 
on $\Z^\infty\oplus \Z^\infty$. 
We have 
\begin{align*}
I^{(2)}-N&=\mat{-I}{-T^+-|S|-X}{-I}{-T^--|S|-X}\\
&=\mat{I}{0}{I}{I}\mat{-I}{0}{0}{T}\mat{I}{T^++|S|+X}{0}{I}
\end{align*}
Since the middle matrix in the above product is injective 
and its cokernel is isomorphic to $G_0$, 
and since the left and right matrices are invertible in $M_2(M_\infty(\Z))$, 
$I^{(2)}-N$ is injective 
and its cokernel is isomorphic to $G_0$. 
In a similar way, 
we can show that $I^{(2)}-M$ is injective 
and its cokernel is isomorphic to $G_1$. 
We define a bijection 
\[
\Z^\infty\oplus\Z^\infty\ni ((x_k),(y_k))\mapsto (z_k)\in \Z^\infty
\]
by $z_{2k-1}=x_k$ and $z_{2k}=y_k$ for $k\in\Z_+$, 
and consider $N,M\in M_\infty(\Z)$ using it. 
We had already seen that $N,M$ satisfy (iv). 
It is easy to see that these satisfy (i) and (ii). 
Since $N_{2k,2k-1},N_{2k-1,2k},N_{2k,2k+2},N_{2k+2,2k}\geq 1$ 
for $k\in\Z_+$, 
we get (iii). 
We are done. 
\end{proof}

\begin{proposition}\label{PropSt}
For countable abelian groups $G_0$ and $G_1$, 
take $N\in M_\infty(\N)$ 
and $M\in M_\infty(\Z)$ satisfying {\em (i)} to {\em (iv)} 
in Lemma \ref{MatrixSt}. 
Then the \Ca $\cO(E_{N,M})$ is a non-unital Kirchberg algebra 
with $K_i(\cO(E_{N,M}))\cong G_i$ for $i=0,1$. 
\end{proposition}

\begin{proof}
First we show that $E_{N,M}$ is minimal. 
Take $k,l\in E^0$ arbitrarily. 
By (iii) in Lemma \ref{MatrixSt}, 
there exists $\mu\in E^*$ 
with $d^*(\mu)=k$ and $r^*(\mu)=l$. 
For $j=1,2,\ldots$, 
we define $\mu_j=((l,l),(l,l),\ldots,(l,l),\mu)$ 
where $(l,l)\in E^1$ is repeated $j$-times. 
Then $\mu_j\in E^*$ satisfies $d^*(\mu_j)=k$ and $r^*(\mu_j)=l$. 
Since $m(l,l)=1$, 
we have $p(\mu_j)=n(l,l)^jp(\mu)$. 
Since $n(l,l)\geq 2$, 
we get 
\[
\sup\{p(\mu)\mid \mu\in E^*, d^*(\mu)=k, r^*(\mu)=l\}
=\infty. 
\]
By Proposition \ref{ET:minimal}, 
$E_{N,M}$ is minimal. 
Since $E$ has a loop, 
$E_{N,M}$ is a Kirchberg algebra 
by Proposition \ref{ET:dico}. 
Since $E^0\times\T$ is not compact, 
$\cO(E_{N,M})$ is not unital. 
Finally we get $K_i(\cO(E))\cong G_i$ 
by Lemma \ref{ENMK} and (iv) in Lemma \ref{MatrixSt}. 
\end{proof}

When a unital \Ca $A$ is in the form $\cO(E\times_{n,m}\T)$, 
Corollary \ref{Cor:unitK} 
implies that both $K_0(\cO(E\times_{n,m}\T))$ and $K_1(\cO(E\times_{n,m}\T))$ 
are finitely generated 
and that the rank of $K_1(\cO(E\times_{n,m}\T))$ 
does not exceed the one of $K_0(\cO(E\times_{n,m}\T))$. 
Thus most of unital Kirchberg algebras do not appear 
by the construction studied so far. 
In Appendix \ref{AppUnital}, 
we see that these two conditions are the only obstructions 
for a unital Kirchberg algebra to be in the form $\cO(E\times_{n,m}\T)$ 
(Proposition \ref{AppProp:unit}). 

We introduce a general way to 
change a given topological graph to a new topological graph 
whose \Ca is simple and unital. 
Combining this and the construction in Proposition \ref{PropSt}, 
we get all unital Kirchberg algebras as \CA s of topological graphs. 

Let $E=(E^0,E^1,d,r)$ be a topological graph 
with non-compact $E^0$. 
Then the \Ca $\cO(E)$ is not unital. 
In \cite[Definition 7.2]{Ka5}, 
we constructed the one-point compactification 
$\tE=(\tE^0,E^1,d,r)$ 
of $E$, 
and showed that $\cO(\tE)$ 
is isomorphic to 
the unitization of $\cO(E)$ 
(\cite[Proposition 7.4]{Ka5}). 
Note that the \Ca $\cO(\tE)$ never be simple. 
We modify this construction to get simple \CA s. 

\begin{definition}
Let $E=(E^0,E^1,d,r)$ be a topological graph 
with non-compact $E^0$. 
For $w\in E^0$, 
we define a topological graph 
$\tE_w=(\tE^0,E^1\amalg \tE^0,d_w,r_w)$ 
such that 
$\tE^0=E^0\cup\{\infty\}$ 
is the one-point compactification of $E^0$,  
and two maps 
$d_w,r_w\colon E^1\amalg \tE^0\to \tE^0$ 
are defined by 
\[
d_w|_{E^1}=d, \quad 
d_w|_{\tE^0}=\id_{\tE^0}, \quad 
r_w|_{E^1}=r, \quad \text{and }\quad r_w(\tE^0)=\{w\}.
\]
\end{definition}

In the statement and the proof of the next proposition, 
note that the positive orbit space $\Orb^+(w)$ 
of $w\in \tE^0$ considered 
in the topological graph $\tE_w$ 
coincides with the one of $w\in E^0$ 
considered in the original topological graph $E$. 

\begin{lemma}\label{Lem:unitiz1}
The topological graph $\tE_w$ is minimal 
if and only if $\Orb^+(w)$ is dense in $E^0$. 
In this case, 
the \Ca $\cO(\tE_w)$ is simple. 
\end{lemma}

\begin{proof}
Suppose that $\Orb^+(w)$ is dense in $E^0$. 
Then for every $v\in \tE^0$, 
$\Orb^+(v)$ is dense in $\tE^0$ 
because $w\in \Orb^+(v)$ and $E^0$ is dense in $\tE^0$. 
Hence $\tE_w$ is minimal. 
Conversely suppose that $\tE_w$ is minimal. 
Since $r_w^{-1}(\infty)=\emptyset$, 
we have $\infty\in (\tE_w)_{\rs{sg}}^0$. 
Hence the minimality of $\tE_w$ implies 
that $\Orb^+(\infty)$ is dense in $\tE^0$. 
Since $\Orb^+(\infty)=\{\infty\}\cup\Orb^+(w)$, 
we see that $\Orb^+(w)$ is dense in $E^0$. 
This completes the proof of the former part. 
Since $\tE^0$ is not discrete, 
$\tE_w$ is not generated by a loop. 
Hence Proposition \ref{PropSimple} implies the latter. 
\end{proof}

\begin{lemma}\label{Lem:unitiz2}
If $E$ is contracting, 
then $\tE_w$ is also contracting. 
\end{lemma}

\begin{proof}
Easy to see. 
\end{proof}

\begin{lemma}\label{Lem:unitiz3}
Suppose $w\in\s{E}{0}{rg}$. 
Let $\varphi_w\colon K_0(C_0(\s{E}{0}{rg}))\to\Z$ 
be the map induced by the \shom 
$C_0(\s{E}{0}{rg})\ni f\mapsto f(w)\in\C$ 
and the natural isomorphism $K_0(\C)\cong\Z$. 
Then there exists an exact sequence 
\[
\begin{CD}
K_0(C_0(\s{E}{0}{rg})) @>>(\iota_*-[\pi_r])\oplus(-\varphi_w)> 
K_0(C_0(E^0)) \oplus \Z 
@>>>  K_0(\cO(\tE_w)) \\
@AAA @. @VVV \\
K_1(\cO(\tE_w)) @<<< K_1(C_0(E^0)) 
@<\iota_*-[\pi_r]<< K_1(C_0(\s{E}{0}{rg})),
\end{CD}
\]
such that $[1]\in K_0(\cO(\tE_w))$ 
is the image of $(0,1)\in K_0(C_0(E^0))\oplus \Z $. 
\end{lemma}

\begin{proof}
When $w\in\s{E}{0}{rg}$, 
we have $(\tE_w)^0_{\rs{rg}}=\s{E}{0}{rg}$. 
Two injections $K_0(C_0(E^0))\to K_0(C(\tE^0))$ 
and $\Z\ni n\mapsto n[1]\in K_0(C(\tE^0))$ 
give an isomorphism 
$K_0(C_0(E^0))\oplus \Z\cong K_0(C(\tE^0))$. 
The group $K_1(C(\tE^0))$ 
is naturally isomorphic to $K_1(C_0(E^0))$. 
Under these isomorphisms, 
the exact sequences in Proposition \ref{Kgroup} 
becomes the desired one.  
\end{proof}

Now we mix the above idea with 
the construction from the previous section 
to get all unital Kirchberg algebras 
as \CA s of topological graphs. 

\begin{lemma}\label{MatrixUn}
For countable abelian groups $G_0, G_1$ and an element $g\in G_0$, 
there exist $N\in M_\infty(\N)$ 
and $M\in M_\infty(\Z)$ satisfying {\rm (i)}, {\rm (ii)}, {\rm (iii)} 
in Lemma \ref{MatrixSt} and  
\begin{enumerate}
\setcounter{enumi}{3}
\renewcommand{\labelenumi}{{\rm (\roman{enumi})'}}
\renewcommand{\itemsep}{0pt}
\item there exist exact sequences 
\[
\begin{CD}
0 @>>> \Z^\infty @>(I-N)\oplus(-\varphi)>> 
\Z^\infty \oplus \Z @>\pi>> G_0 @>>> 0
\end{CD}
\]
\[
\begin{CD}
0 @>>> \Z^\infty @>I-M>> \Z^\infty @>>> G_1 @>>> 0
\end{CD}
\]
such that 
$\varphi\colon\Z^\infty \ni (x_k)_{k=1}^\infty\mapsto x_1\in \Z$ 
and $\pi(0,1)=g$. 
\end{enumerate}
\end{lemma}

\begin{proof}
Take a surjective homomorphism 
$\pi_0'\colon \Z^\infty\to G_0$ 
such that $\ker \pi_0'$ has an infinite rank. 
Fix an isomorphism $T'$ from 
$\{(x_k)_{k=1}^\infty\in \Z^\infty\mid x_1=0\}$ 
to $\ker \pi_0'$. 
Take $a\in \Z^\infty$ with $\pi_0'(a)=g$, 
and define $T\colon\Z^\infty\to\Z^\infty$ by 
\[
T((x_k)_{k=1}^\infty)=-x_1a+T'((0,x_2,x_3,\ldots)). 
\]
We define $\pi_0\colon \Z^\infty\oplus\Z\to G_0$ 
by $\pi_0(x,n)=\pi_0'(x)+ng$. 
Then the sequence 
\[
\begin{CD}
0 @>>> \Z^\infty @>T\oplus \varphi>> 
\Z^\infty \oplus \Z@>\pi_0>> G_0 @>>> 0
\end{CD}
\]
is exact and the image of $(0,1)$ is $g\in G_0$. 
Take $S\in M_\infty(\Z)$ 
which is injective and whose cokernel is isomorphic to $G_1$. 
Using these $T,S$, we get $N\in M_\infty(\N)$ 
and $M\in M_\infty(\Z)$ by the same way as 
in the proof of Lemma \ref{MatrixSt}. 
These $N,M$ satisfy (i), (ii), (iii) and (iv)'. 
\end{proof}

For countable abelian groups $G_0, G_1$ and an element $g\in G_0$, 
take $N\in M_\infty(\N)$ 
and $M\in M_\infty(\Z)$ satisfying the conditions 
in Lemma \ref{MatrixUn}. 
We set $F=E_{N,M}$ and $w=(1,1)\in \Z_+\times\T=F^0$. 

\begin{proposition}\label{PropUn}
The \Ca $\cO(\widetilde{F}_w)$ is a unital Kirchberg algebra with 
\[
\big(K_0(\cO(\widetilde{F}_w)),[1],K_1(\cO(\widetilde{F}_w))\big)
\cong \big(G_0,g,G_1\big).
\]
\end{proposition}

\begin{proof}
By using Lemma \ref{MatrixUn}, 
we see that the topological graph $F$ is minimal and contracting 
in a similar way to the proof of Proposition \ref{PropSt}. 
Hence $\cO(\widetilde{F}_w)$ is a unital Kirchberg algebra 
by Lemmas \ref{Lem:unitiz1} and \ref{Lem:unitiz2}. 
Now we get 
\[
\big(K_0(\cO(\widetilde{F}_w)),[1],K_1(\cO(\widetilde{F}_w))\big)
\cong \big(G_0,g,G_1\big)
\]
by Lemma \ref{Lem:unitiz3} and the condition (iv)' in Lemma \ref{MatrixUn}. 
\end{proof}

By Propositions \ref{PropSt} and \ref{PropUn}, 
we get Theorem C. 
In \cite{Ka9}, 
we see that the construction here can be used to produce 
various actions on Kirchberg algebras.

\appendix

\section{Unital Kirchberg algebras}\label{AppUnital}

As explained in the remark after Proposition \ref{PropSt}, 
if a unital Kirchberg algebra $A$ is in the form $\cO(E\times_{n,m}\T)$, 
then both $K_0(A)$ and $K_1(A)$ 
are finitely generated abelian groups, 
and the rank of $K_1(A)$ 
does not exceed the one of $K_0(A)$. 
We will prove the converse. 

Take finitely generated abelian groups $G_0$ and $G_1$ 
such that the rank of $G_1$ 
does not exceed the one of $G_0$. 
Let $l_0\in\N$ be the difference of the ranks of $G_0$ and $G_1$. 
Then, there exist $L\in\Z_+$, 
$T\in M_{L,L-l_0}(\N)$ and $S\in M_L(\N)$ with 
\[
G_0\cong \coker T \oplus \ker S,\quad 
G_1\cong \ker T \oplus \coker S.
\]
Take $\widetilde{T}\in M_L(\N)$ 
such that its restriction to $M_{L,L-l_0}(\N)$ is $T$ 
(for example, set $\widetilde{T}_{k,l}=0$ for $l>L-l_0$). 
We set $X\in M_L(\N)$ by 
\[
X_{k,l}=\begin{cases}
1 & \text{if }|k-l|\leq 1,\\
0 & \text{if }|k-l|\geq 2.
\end{cases}
\]
Let us denote by $I$ the identity matrix in $M_{L}(\N)$. 
We set $N, M\in M_{2L}(\N)$ by 
\[
\widetilde{N}=\mat{2I}{\widetilde{T}+S+X}{I}{I+S+X},\quad
M=\mat{I}{S}{I}{I}.
\]
We denote by $N$ the restriction of $\widetilde{N}$ to $M_{2L,2L-l_0}(\N)$. 
Let $I_1$ be the identity matrix in $M_{2L}(\N)$, 
and $I_0$ be its restriction to $M_{2L,2L-l_0}(\N)$. 
In a similar way to the proof of Lemma \ref{MatrixSt}, 
we see 
\begin{align*}
\coker (I_0-N)&\cong \coker T,& 
\ker (I_0-N)&\cong \ker T,\\ 
\coker (I_1-M)&\cong \coker S,& 
\ker (I_1-M)&\cong \ker S. 
\end{align*}
We define a discrete graph $E=(E^0,E^1,d,r)$ as follows. 
Set $E^0=\{1,2,\ldots,2L\}$. 
Set $\Omega=\{(k,l)\in E^0\times E^0\mid \widetilde{N}_{k,l}\geq 1\}$, 
and 
\[
E^1=\Omega\amalg\{e_{l,p}\mid 2L-l_0+1\leq l\leq 2L,\ p=1,2,\ldots\}. 
\]
We define $d,r\colon E^1\to E^0$ by $d(k,l)=k$, $r(k,l)=l$ 
for $(k,l)\in \Omega$ and $d(e_{l,p})=1$, $r(e_{l,p})=l$. 
We define $n,m\colon E^1\to \N$ 
by $n(k,l)=\widetilde{N}_{k,l}$, $m(k,l)=M_{k,l}$ 
for $(k,l)\in \Omega$ and $n(e_{l,p})=1$, $m(e_{l,p})=0$. 
In a very similar way to the proof of Proposition \ref{PropSt}, 
we see that $\cO(E\times_{n,m}\T)$ is a Kirchberg algebra. 
We compute its $K$-groups. 
We have $\s{E}{0}{rg}=\{1,2,\ldots,2L-l_0\}$ 
and $E_m^0=\{1,2,\ldots,2L\}$. 
The matrices $I_0,I_1,N,M$ 
defined before Proposition \ref{Prop:K} 
are nothing but the ones defined here. 
Hence Proposition \ref{Prop:K} and the above computation show that 
the unital Kirchberg algebra $A=\cO(E\times_{n,m}\T)$ satisfies 
$K_i(A)\cong G_i$ for $i=0,1$. 
Since $(1,\ldots,1,1,\ldots,1)$ is 
the image of $(-1,\ldots,-1,0,\ldots,0)\in\Z^{2L-l_0}$ 
(here $-1$ is repeated $L$-times, and $0$ is repeated $(L-l_0)$-times) 
under $I_0-N$, 
we have $[1_{A}]=0$. 

Now take $G_0$ and $G_1$ as above, 
and take $g\in G_0$ with $g\neq 0$. 
We will construct a discrete graph $E'$ and the two maps $n',m'$ 
using $\widetilde{N},M$ defined above 
so that the \Ca $A'=\cO(E'\times_{n',m'}\T)$ 
is a unital Kirchberg algebra with 
\[
\big(K_0(A'),[1_{A'}],K_1(A')\big)
\cong (G_0,g,G_1).
\]
Take $\widetilde{N},M\in M_{2L}(\N)$ as above. 
We may assume that $\ker (I_1-M)=0$ 
by choosing $T,S$ in the beginning of the above proof 
so that $\ker S=0$. 
Hence there exists a surjection $\pi\colon \Z^{2L}\to G_0$ 
such that $\ker\pi$ is the image of $I_0-N$. 
As shown above, $\pi(1,1,\ldots,1)=0$. 
Take $a=(a_1,a_2,\ldots,a_{2L})\in \Z^{2L}$ with $\pi(a)=g$. 
By adding $(b,b,\ldots,b)\in \ker\pi\subset \Z^{2L}$ 
for a suitable $b\in\Z$, 
we may assume $\min\{a_k\mid k=1,2,\ldots,2L\}=0$. 
Then there exists $k$ with $a_k\geq 1$ because $g\neq 0$. 
Take $k_0$ with $a_{k_0}=0$. 
We define 
$\widetilde{N}'=(\widetilde{N}'_{k,l})_{k,l=0}^{2L}\in M_{1+2L}(\N)$ 
by
\[
\widetilde{N}'_{k,l}=\begin{cases}
2&\text{for $k=0$ and $l=0$,}\\
2\widetilde{N}_{k_0,k_0}-2&\text{for $k=0$ and $l=k_0$,}\\
2\widetilde{N}_{k_0,l}&\text{for $k=0$ and $l\neq 0,k_0$,}\\
a_k&\text{for $k\neq 0$ and $l=0$,}\\
\widetilde{N}_{k,l}&\text{for $1\leq k,l\leq 2L$.}
\end{cases}
\]
Let $N'$ be the restriction of $\widetilde{N}'$ 
to $M_{1+2L,1+2L-l_0}(\N)$. 
Let $I_1'$ be the identity matrix in $M_{1+2L}(\N)$, 
and $I_0'$ be its restriction to $M_{1+2L,1+2L-l_0}(\N)$. 
We define 
\[
\pi'\colon \Z^{1+2L}\ni (x_0,x_1,\ldots,x_{2L})\mapsto 
\pi(x_1,\ldots,x_{2L})+(2x_{k_0}-x_0)g \in G. 
\]
Then one can verify that $\pi'$ is a surjective homomorphism, 
$\pi'(1,1,\ldots,1)=g$, 
and $\ker\pi'$ is the image of $I_0'-N'\colon \Z^{1+2L-l_0}\to \Z^{1+2L}$. 
Thus $\pi'$ defines an isomorphism 
$\coker (I_0'-N')\cong G_0$ which sends the image of 
$(1,1,\ldots,1)\in \Z^{1+2L}$ to $g\in G$. 
Similarly, 
we define $M'\in M_{1+2L}(\N)$ by 
\[
M'_{k,l}=\begin{cases}
2&\text{for $k=0$ and $l=0$,}\\
0&\text{for $k=0$ and $l\neq 0$,}\\
1&\text{for $k\neq 0$ and $l=0$,}\\
M_{k,l}&\text{for $1\leq k,l\leq 2L$.}
\end{cases}
\]
Then $\ker (I_1'-M')\cong\ker (I_1-M)=0$, 
and we have an isomorphism 
$G_1\cong \ker (I_0'-N') \oplus \coker (I_1'-M')$.
From the two matrices $\widetilde{N}',M'\in M_{1+2L}(\N)$, 
we define a discrete graph $E'$ and the two maps $n',m'$ 
similarly as above. 
In a similar way to the proof of Proposition \ref{PropSt}, 
one can show 
that the \Ca $A'=\cO(E'\times_{n',m'}\T)$ 
is a unital Kirchberg algebra satisfying 
\[
\big(K_0(A'),[1_{A'}],K_1(A')\big)
\cong (G_0,g,G_1)
\]
(because $\widetilde{N}'_{0,0}=M'_{0,0}=2$ 
we need to modify the proof). 
Now we have shown the following. 

\begin{proposition}\label{AppProp:unit}
A unital Kirchberg algebra $A$ is 
in the form $\cO(E\times_{n,m}\T)$ 
if and only if both $K_0(A)$ and $K_1(A)$ 
are finitely generated, 
and the rank of $K_1(A)$ 
does not exceed the one of $K_0(A)$. 
\end{proposition}

By the above construction, 
we prove the following. 
Recall that a discrete graph $E=(E^0,E^1,d,r)$ 
is said to be {\em finite} 
if both $E^0$ and $E^1$ are finite sets. 

\begin{proposition}\label{AppProp:finite}
A unital Kirchberg algebra $A$ is 
in the form $\cO(E\times_{n,m}\T)$ 
for a finite discrete graph $E$ 
if and only if $K_0(A)$ and $K_1(A)$ 
are finitely generated abelian groups with same rank. 
\end{proposition}

\section{Generators of $\cO(E\times_{n,m}\T)$ and their relations}\label{app3}

Recall that for a topological graph $E$ 
a pair $(T^0,T^1)$ of maps 
satisfying the three conditions in Definition \ref{DefO(E)} 
is called a Cuntz-Krieger $E$-pair. 
When a pair $(T^0,T^1)$ only satisfies (i) and (ii), 
then we call it a Toeplitz $E$-pair. 
The \Ca generated by the universal Toeplitz $E$-pair 
is denoted by $\cT(E)$. 

Take a discrete graph $E=(E^0,E^1,d,r)$ 
and two maps $n\colon E^1\to\Z_+$ and $m\colon E^1\to\Z$. 
In this appendix, 
we will present a generator of the \Ca $\cO(E\times_{n,m}\T)$ 
and its relations. 
For each $v\in E^0$, 
let $p_v,u_v\in \cO(E\times_{n,m}\T)$ 
be the images of the unit 
and the generating unitary of 
$C(\{v\}\times\T)\subset C_0(E^0\times\T)$ 
under the \shom $t^0\colon C_0(E^0\times\T)\to \cO(E\times_{n,m}\T)$. 
For $e\in E^1$ and $k\in\Z$, 
we define $s_{e,k}=t^1(\xi_{e,k})$ 
where $\xi_{e,k}\in C_\td(E^1\times\T)$ 
is defined by $\xi_{e,k}(e',z)=0$ for $e'\neq e$ and 
$\xi_{e,k}(e,z)=z^k/\sqrt{n(e)}$.

\begin{lemma}\label{Lem:relation}
The set of elements $\{p_v,u_v\}_{v\in E^0}$ 
and $\{s_{e,k}\}_{e\in E^1,k\in\Z}$ in $\cO(E\times_{n,m}\T)$ 
satisfies the following. 
\begin{enumerate}
\rom
\item $u_v^*u_v=u_vu_v^*=p_v$ for $v\in E^0$, 
\item $\{p_v\}_{v\in E^0}$ are mutually orthogonal projections, 
\item $s_{e,k}^*s_{e,k}=p_{d(e)}$ for $e\in E^1,k\in\Z$, 
\item $\{s_{e,k}s_{e,k}^*\}_{e\in E^1,0\leq k<n(e)}$ 
are mutually orthogonal projections, 
\item $s_{e,k}u_{d(e)}=s_{e,k+n(e)}$ for $e\in E^1,k\in\Z$, 
\item $u_{r(e)}s_{e,k}=s_{e,k+m(e)}$ for $e\in E^1,k\in\Z$, 
\item $p_v=\sum_{e\in r^{-1}(v)}\sum_{k=0}^{n(e)-1}s_{e,k}s_{e,k}^*$ 
for $v\in \s{E}{0}{rg}\cap E^0_m$, 
\item $p_v-u_v=\sum_{e\in r^{-1}(v)\setminus m^{-1}(0)}
\sum_{k=0}^{n(e)-1}(s_{e,k}-s_{e,k+m(e)})s_{e,k}^*$ 
for $v\in \s{E}{0}{sg}\cap E^0_m$. 
\end{enumerate}
\end{lemma}

\begin{proof}
It follows from the direct computation. 
\end{proof}

A partial isometry whose initial and final projections coincide 
is called a {\em partial unitary}. 
For a partial unitary $u$ with $u^*u=uu^*=p$, 
we set $u^0=p$ and $u^{-n}=(u^*)^{n}$ for $n\geq 1$. 

The \Ca $\cO(E\times_{n,m}\T)$ is the universal \Ca 
generated by $\{p_v,u_v\}_{v\in E^0}$ 
and $\{s_{e,k}\}_{e\in E^1,k\in\Z}$ 
whose relations are listed in Lemma \ref{Lem:relation}. 
By reducing the number of generators, 
we get the following. 

\begin{proposition}\label{Prop:univ}
Let $B$ be a \Ca generated 
by a family $\{U_v\}_{v\in E^0}$ of partial unitaries with 
orthogonal ranges 
and a family $\{S_{e,k}\}_{e\in E^1,0\leq k<n(e)}$ 
of partial isometries with 
orthogonal ranges 
satisfying the following; 
\begin{enumerate}
\rom
\item $S_{e,k}^*S_{e,k}=U_{d(e)}^*U_{d(e)}$ for $e\in E^1$ and $0\leq k<n(e)$, 
\item $U_{r(e)}S_{e,k}=S_{e,k+m(e)}$ 
for $e\in E^1$ and $0\leq k<n(e)$, 
\item $U_v^*U_v=\sum_{e\in r^{-1}(v)}\sum_{k=0}^{n(e)-1}S_{e,k}S_{e,k}^*$ 
for $v\in \s{E}{0}{rg}\cap E^0_m$, 
\item $U_v^*U_v-U_v=\sum_{e\in r^{-1}(v)\setminus m^{-1}(0)}
\sum_{k=0}^{n(e)-1}(S_{e,k}-S_{e,k+m(e)})S_{e,k}^*$ 
for $v\in \s{E}{0}{sg}\cap E^0_m$, 
\end{enumerate}
where in {\rm (ii)} and {\rm (iv)} 
$S_{e,k+m(e)}$ is defined by 
$S_{e,k+m(e)}=S_{e,k'}U_{d(e)}^{l}$ 
using unique $k'\in \{0,1,\ldots,n(e)-1\}$ and $l\in\Z$ 
with $k+m(e)=k'+n(e)l$. 
Then there exists a \shom 
$\rho\colon \cO(E\times_{n,m}\T)\to B$ 
with $\rho(u_v)=U_v$ and $\rho(s_{e,k})=S_{e,k}$. 
\end{proposition}

\begin{proof}
We can define a \shom 
$T^0\colon C_0(E^0\times\T)\to B$ 
by sending the generating unitary of 
$C(\{v\}\times\T)\subset C_0(E^0\times\T)$ 
to $U_v$ for each $v\in E^0$. 
For $k\in \{0,1,\ldots,n(e)-1\}$ and $l\in\Z$, 
we set $S_{e,k+n(e)l}=S_{e,k}U_{d(e)}^{l}$.
We define a linear map $T^1\colon C_{\td}(E^1\times\T)\to B$ 
by sending $\xi_{e,k}\in C_\td(E^1\times\T)$ 
defined by $\xi_{e,k}(e',z)=0$ for $e'\neq e$ and 
$\xi_{e,k}(e,z)=z^k/\sqrt{n(e)}$ 
to $S_{e,k}$ for $e\in E^1$ and $k\in\Z$. 
By (i), 
the map $T^1$ is well defined 
and satisfies $T^1(\xi)^*T^1(\eta)=T^0(\ip{\xi}{\eta})$ 
for $\xi,\eta\in C_\td(E^1\times\T)$. 
By (ii), 
the pair $(T^0,T^1)$ is a Toeplitz $(E\times_{n,m}\T)$-pair, 
and one can verify that 
it is a Cuntz-Krieger $(E\times_{n,m}\T)$-pair 
from (iii) and (iv). 
Thus we get a \shom 
$\rho\colon \cO(E\times_{n,m}\T)\to B$ 
with $\rho(u_v)=U_v$ and $\rho(s_{e,k})=S_{e,k}$. 
\end{proof}

\begin{remark}\label{Rem:univ}
By the above proof, 
we see that 
if the families $\{U_v\}$ and $\{S_{e,k}\}$ only satisfy (i) and (ii), 
then we get a \shom $\rho\colon \cT(E\times_{n,m}\T)\to B$. 
\end{remark}

\begin{corollary}\label{AppCor:finpre}
For a finite graph $E$, 
the \Ca $\cO(E\times_{n,m}\T)$ 
is finitely presented, 
that is, 
$\cO(E\times_{n,m}\T)$ is the universal \Ca 
generated by finite elements satisfying finite relations. 
\end{corollary}

Note that 
we have $\s{E}{0}{sg}\cap E^0_m=\emptyset$ for a finite graph $E$. 
Hence in this case 
the condition (iv) in Proposition \ref{Prop:univ} 
is void. 

\begin{proposition}
A Kirchberg algebra $A$ such that 
$K_0(A)$ and $K_1(A)$ are finitely generated abelian groups 
with same rank is finitely presented. 
\end{proposition}

\begin{proof}
Combine Proposition \ref{AppProp:finite} 
and Corollary \ref{AppCor:finpre}. 
\end{proof}

\begin{example}\label{ExUniv1}
Let $A$ be the unital Kirchberg algebra with 
\[
\big(K_0(A),[1_{A}],K_1(A)\big)
\cong (\Z/p\Z,1,\Z/q\Z)
\]
for two integers $p,q$ with $p,q\geq 2$. 
Set $n=1+p\in \Z_+$. 
Then either $1+q$ or $1-q$ is not in $n\Z$. 
Let $m$ be such an integer. 
Then the \Ca $\cO(E_{n,m})$ considered in Section \ref{SecEnm} 
is isomorphic to $A$ by Proposition \ref{KOEnm} 
and Proposition \ref{EnmMin}. 
Thus by Proposition \ref{Prop:univ} 
the \Ca $A$ is the universal \Ca generated by 
a unitary $u$ and a family $\{s_k\}_{k=0}^{n-1}$ 
of isometries such that $\sum_{k=0}^{n-1}s_ks_k^*=1$ 
and $us_k=s_{k'}u^l$ 
for $k'\in \{0,1,\ldots,n-1\}$ and $l\in\Z$ 
with $k+m=k'+nl$. 
If we set $s_{k+nl}=s_ku^l$ and set $s_{i,k}=s_{i+mk}$ 
for $i=0,1,\ldots,d-1$ and $k=0,1,\ldots,n/d-1$ 
where $d=(n,|m|)$, 
then we have $us_{i,k}=s_{i,k+1}$ for $0\leq k<n/d-1$ 
and $us_{i,n/d-1}=s_{i,0}u^{m/d}$. 
Thus the \Ca $A$ is also the universal \Ca generated by 
a unitary $u$ and a family $\{s_{i,k}\}_{0\leq i<d, 0\leq k<n/d}$ 
of isometries satisfying $\sum_{i,k}s_{i,k}s_{i,k}^*=1$ 
and the two relations above. 
In particular, 
when $n$ and $|m|$ are mutually prime, 
the \Ca $A$ is the universal \Ca generated by 
two elements $u$ and $s$ with 
\begin{enumerate}
\rom
\item $u^*u=uu^*=s^*s=\sum_{k=0}^{n-1}u^kss^*(u^*)^{k}=1$, 
\item $u^ns=su^m$. 
\end{enumerate}
\end{example}

\begin{example}
Let $A$ be the unital Kirchberg algebra with 
\[
\big(K_0(A),[1_{A}],K_1(A)\big)
\cong (\Z^n,0,\Z^n\oplus F)
\]
where $n\in\N$ and $F$ is a finite abelian group. 
Take natural numbers $q_1,\ldots,q_K$ 
such that $K_1(A)\cong\bigoplus_{k=1}^K \Z/q_k\Z$. 
By adding $1$ to the list $q_1,\ldots,q_K$ if necessary, 
we may assume that $K$ is even. 
Then the two matrices $N,M\in M_K(\N)$ defined as 
\[
N=\left(\begin{array}{ccccc}
3 & 1 & & & \\
  & 2 & 1 & & \\
 & & 2 & \ddots &\phantom{\ddots} \\
 &\phantom{\ddots} & & \ddots & 1 \\
1 & & & & 2
\end{array}
\right),\quad
M=\left(\begin{array}{ccccc}
1 & q_1 & & & \\
  & 1 & q_2 & & \\
 & & 1 & \ddots & \\
 & & & \ddots & q_{K-1} \\
q_K & & & & 1
\end{array}
\right)
\]
satisfy the analogous conditions of (i), (ii), (iii) in Lemma \ref{MatrixSt} 
and  
\[
\ker (I-N)=\coker (I-N)=0,\ \ker (I-M)\cong K_0(A), 
\text{ and }\coker (I-M)\cong K_1(A).
\] 
Then the \Ca $\cO(E_{N,M})$ is isomorphic to $A$ 
(see Remark \ref{RemOENM} and the proof of Proposition \ref{PropSt}). 
By Proposition \ref{Prop:univ}, 
the \Ca $A$ is the universal \Ca generated by 
by a family $\{u_k\}_{k=1}^K$ of partial unitaries with 
orthogonal ranges and 
a family $\{s_k,t_k\}_{k=1}^K$ of partial isometries 
satisfying the following; 
\begin{enumerate}
\rom
\item $s_k^*s_k=t_k^*t_k=u_k^*u_k$ for $1\leq k\leq K$, 
\item $u_1^3s_1=s_1u_1$, $u_1t_K=t_Ku_K^{q_K}$, 
and $u_k^2s_k=s_ku_k$, 
$u_kt_{k-1}=t_{k-1}u_{k-1}^{q_{k-1}}$ for $2\leq k\leq K$, 
\item $u_1u_1^*=s_1s_1^*+u_1s_1s_1^*u_1^*+u_1^2s_1s_1^*(u_1^*)^2+t_Kt_K^*$ and 
$u_ku_k^*=s_ks_k^*+u_ks_ks_k^*u_k^*+t_{k-1}t_{k-1}^*$ for $2\leq k\leq K$. 
\end{enumerate}
\end{example}

\begin{remark}\label{remsemiproj}
The finite presentations above may be used to prove the following conjecture 
(for details, see \cite{B2}); 

\medskip
\noindent 
{\em Conjecture} (\cite[Conjecture 3.6]{B2}). 
A Kirchberg algebra is semiprojective if and only if 
its $K$-theory is finitely generated. 
\medskip

As Blackadar explained in \cite{B2}, 
one approach to the conjecture 
is finding a finite presentation of a Kirchberg algebra 
whose $K$-theory is finitely generated, 
and showing that the relations are stable. 
The relation (i) in Example \ref{ExUniv1} is stable, 
but the relation (ii) is not on its own. 
However, 
it seems to be reasonable to conjecture that 
the pair of the two relations is stable. 
In fact, the pair is stable when $m=1$ 
because the universal \Ca $\cO(E_{n,m})$ in this case 
is the unital Kirchberg algebra $A$ with 
\[
\big(K_0(A)),[1_{A}],K_1(A)\big)
\cong (\Z\oplus \Z/(n-1)\Z,(0,1),\Z)
\]
by Proposition \ref{KOEnm} 
and Proposition \ref{EnmMin}, 
which is isomorphic to a Cuntz-Krieger algebra \cite{CK}, 
hence is semiprojective. 
\end{remark}

\section{A proof of Lemma \ref{Lem:appendix}}\label{app1}

In this appendix, 
we will prove Lemma \ref{Lem:appendix}. 
Let $E=(E^0,E^1,d,r)$ be a discrete graph and 
$n\colon E^1\to\Z_+$, $m\colon E^1\to\Z$ be two maps. 
Let $F^0\subset E^0$ and $F^1\subset E^1$ be 
finite subsets satisfying $d(F^1),r(F^1)\subset F^0$. 
The restrictions of $n,m\colon E^1\to\Z$ to $F^1\subset E^1$ 
are also denoted by the same symbols $n,m$. 
We get a topological graph $F\times_{n,m}\T$. 
One can easily check that 
the pair $T=(T^0,T^1)$ of the restrictions of $t^0,t^1$ 
to $C(F^0\times\T)\subset C_0(E^0\times\T)$ 
and $C(F^1\times\T)\subset C_{\td}(E^1\times\T)$ respectively, 
is a Toeplitz $(F\times_{n,m}\T)$-pair. 
Hence we get a surjection from $\cT(F\times_{n,m}\T)$ 
to the \Csa of $\cO(E\times_{n,m}\T)$ 
generated by $t^0(C(F^0\times\T))$ and $t^1(C(F^1\times\T))$. 
Therefore Lemma \ref{Lem:appendix} follows from the next lemma. 

\begin{lemma}\label{Lem:forappendix} 
For a finite graph $F=(F^0,F^1,d,r)$ with no loops 
and two maps $n\colon F^1\to\Z_+$, $m\colon F^1\to\Z$, 
the \Ca $\cT(F\times_{n,m}\T)$ is a circle algebra. 
\end{lemma}

We will prove Lemma \ref{Lem:forappendix}. 
Take a finite graph $F=(F^0,F^1,d,r)$ with no loops 
and two maps $n\colon F^1\to\Z_+$, $m\colon F^1\to\Z$. 
Note that the path space $F^*$ of $F$ is a finite set. 
Define $n\colon F^*\to\Z_+$ by $n(v)=1$ for $v\in F^0$, 
and $n(\mu)=n(e_1)\cdots n(e_k)$ 
for $\mu=(e_1,\ldots,e_k)\in F^k$. 
We define a finite set $\lambda$ by 
\[
\lambda=\big\{(\mu,k)\ \big|\ \mu\in F^*,\ 
k\in \{0,1,\ldots,n(\mu)-1\}\big\}. 
\]
Let $\{w_{(\mu,k),(\nu,l)}\}_{(\mu,k),(\nu,l)\in \lambda}$ 
be the matrix units of $\M_{|\lambda|}$, 
and $A$ be the \Csa of $\M_{|\lambda|}$ 
spanned by 
$\{w_{(\mu,k),(\nu,l)}\mid d^*(\mu)=d^*(\nu)\}$. 
We denote by $u$ the generating unitary of 
$1_{A}\otimes C(\T)\subset A\otimes C(\T)$. 
We consider $A$ as a unital \Csa of $A\otimes C(\T)$. 
Thus $u$ commutes all elements in $A$. 
We will show that $\cT(F\times_{n,m}\T)$ 
is isomorphic to $A\otimes C(\T)$. 

For $\mu\in F^*$, 
we define $P_{\mu}\in A\otimes C(\T)$ by 
\[
P_{\mu}=\sum_{l=0}^{n(\mu)-1}w_{(\mu,l),(\mu,l)}. 
\]
For $e\in F^1$ and $\mu\in F^*$ with $d(e)=r^*(\mu)$ 
and $k\in\{0,1,\ldots,n(e)-1\}$, 
we define $S_{(e,\mu),k}\in A\otimes C(\T)$ by 
\[
S_{(e,\mu),k}=\sum_{l=0}^{n(\mu)-1}w_{((e,\mu),kn(\mu)+l),(\mu,l)}. 
\]

It is routine to see the following. 

\begin{lemma}
The family $\{P_{\mu}\}_{\mu\in F^*}$ 
is a family of mutually orthogonal projections, 
and for all $e\in F^1$ and $\mu\in F^*$ with $d(e)=r^*(\mu)$, 
we have $S_{(e,\mu),k}^*S_{(e,\mu),k}=P_{\mu}$ 
for $k\in \{0,1,\ldots,n(e)-1\}$ and 
$\sum_{k=0}^{n(e)-1}
S_{(e,\mu),k}S_{(e,\mu),k}^*=P_{(e,\mu)}$. 
\end{lemma}

We define a partial unitary $U_{\mu}\in A\otimes C(\T)$ with 
$U_{\mu}^*U_{\mu}=U_{\mu}U_{\mu}^*=P_{\mu}$ 
for $\mu\in F^*$ recursively by 
$U_{\mu}=w_{(v,0),(v,0)}u$ for $\mu=v\in F^0$, and 
\[
U_{(e,\mu)}=\Big(
S_{(e,\mu),0}U_{\mu}{S_{(e,\mu),n(e)-1}}^*
+ \sum_{k=0}^{n(e)-2}S_{(e,\mu),k+1}{S_{(e,\mu),k}}^*
\Big)^{m(e)}
\]
for $e\in F^1$ and $\mu\in F^*$ with $d(e)=r^*(\mu)$. 
If we define $S_{(e,\mu),k+n(e)l}=S_{(e,\mu),k}U_{\mu}^l$ 
for $k\in \{0,1,\ldots,n(e)-1\}$ and $l\in\Z$, 
then one can check that 
\[
U_{(e,\mu)}=\Big(
\sum_{k=0}^{n(e)-1}S_{(e,\mu),k+1}{S_{(e,\mu),k}}^*
\Big)^{m(e)}
=\sum_{k=0}^{n(e)-1}S_{(e,\mu),k+m(e)}{S_{(e,\mu),k}}^*.
\]

For $v\in F^0$, 
we define $p_v=\sum_{\mu\in F^*, r^*(\mu)=v}P_\mu$ 
and $u_v=\sum_{\mu\in F^*, r^*(\mu)=v}U_\mu$. 
Then $\{p_v\}_{v\in F^0}$ is a family of 
mutually orthogonal projections, 
and $u_v^*u_v=u_vu_v^*=p_v$ for $v\in F^0$. 
For $e\in F^1$ and $k\in \{0,1,\ldots,n(e)-1\}$, 
we define 
\[
s_{e,k}=\sum_{\mu\in F^*, r^*(\mu)=d(e)}S_{(e,\mu),k}. 
\]
Then $\{s_{e,k}s_{e,k}^*\}_{e\in F^1,0\leq k<n(e)}$ is a family of 
mutually orthogonal projections, 
and $s_{e,k}^*s_{e,k}=p_{d(e)}$ 
for $e\in F^1$ and $k\in \{0,1,\ldots,n(e)-1\}$. 
For $k\in \{0,1,\ldots,n(e)-1\}$ and $l\in\Z$, 
we set $s_{e,k+n(e)l}=s_{e,k}u_{d(e)}^l$. 
Then we have 
\[
s_{e,l}=\sum_{\mu\in F^*, r^*(\mu)=d(e)}S_{(e,\mu),l}. 
\]
for $e\in F^1$ and $l\in\Z$. 
Using this equality, 
we can show that $u_{r(e)}s_{e,k}=s_{e,k+m(e)}$ 
for $e\in E^1$ and $k\in \{0,1,\ldots,n(e)-1\}$. 

By the proof of Proposition \ref{Prop:univ}, 
we get a Toeplitz $(F\times_{n,m}\T)$-pair $(T^0,T^1)$ on $A\otimes C(\T)$
where a \shom $T^0\colon C(F^0\times\T)\to A\otimes C(\T)$ 
sends the unit and the generating unitary of 
$C(\{v\}\times\T)$ to $p_v$ and $u_v$, 
and a linear map $T^0\colon C_\td(F^1\times\T)\to A\otimes C(\T)$ 
sends $\xi_{e,k}$ defined by $\xi_{e,k}(e',z)=0$ for $e'\neq e$ and 
$\xi_{e,k}(e,z)=z^k/\sqrt{n(e)}$ 
to $s_{e,k}$. 
This Toeplitz pair $(T^0,T^1)$ induces 
a \shom $\rho\colon\cT(F\times_{n,m}\T)\to A\otimes C(\T)$ 
(see Remark \ref{Rem:univ}). 
We will show that the \shom $\rho$ is an isomorphism. 
To show that $\rho$ is injective, 
it suffices to check the following 
two conditions (see \cite[Corollary 3.22]{Ka5}): 
\begin{itemize}
\item $\rho$ admits a gauge action, and 
\item no non-zero element $f$ of $C_0(F^0\times\T)$ satisfies 
$T^0(f)\neq \varPhi(\pi_{\tr}(f))$, 
\end{itemize}
where $\varPhi\colon \cK(C_{\td}(F^1\times\T))\to A\otimes C(\T)$ 
is the \shom defined by $\varPhi(\theta_{\xi,\eta})=T^1(\xi)T^1(\eta)^*$. 
One can check that the action $\beta\colon\T\curvearrowright A\otimes C(\T)$, 
defined by $\beta_z(w_{(\mu,k),(\nu,l)})=z^{p-q}w_{(\mu,k),(\nu,l)}$ 
for $(\mu,k),(\nu,l)\in \lambda$ with $\mu\in F^p$ and $\nu\in F^q$, 
gives a gauge action for the \shom $\rho$. 
The second condition follows from the next lemma 
because the spectrum of $U_v$ is $\T$ for all $v\in F^0$. 

\begin{lemma}\label{Lem:Uv}
Let us take $v\in F^0$, and 
let $z$ be the generating unitary of 
$C(\{v\}\times\T)$. 
Then we have $T^0(z)=u_v=\sum_{\mu\in F^*, r^*(\mu)=v}U_\mu$, and 
\[
\varPhi(\pi_{\tr}(z))
=\sum_{e\in r^{-1}(v), 0\leq k<n(e)}s_{e,k+m(e)}s_{e,k}^*
=\sum_{\mu\in F^*\setminus F^0, r^*(\mu)=v}U_\mu.
\]
Thus $T^0(z)-\varPhi(\pi_{\tr}(z))=U_v$. 
\end{lemma}

\begin{proof}
Straightforward. 
\end{proof}

Thus the \shom $\rho$ is injective. 
We will show that it is surjective. 
Take $(\mu,k),(\nu,l)\in\lambda$ with $d^*(\mu)=d^*(\nu)$. 
Take $e\in F^1$ with $d(e)=r^*(\mu)$ 
and $k'\in \{0,1,\ldots,n(e)-1\}$. 
Then we have 
\[
w_{((e,\mu),k'n(e)+k),(\nu,l)}u
=S_{(e,\mu),k'}(w_{(\mu,k),(\nu,l)}u)
=s_{e,k'}(w_{(\mu,k),(\nu,l)}u). 
\]
By Lemma \ref{Lem:Uv}, 
$w_{(v,0),(v,0)}u=U_v$ is in the image of $\rho$ 
for all $v\in F^0$. 
Combining the two facts above, 
we can show by induction 
that $w_{(\mu,k),(\nu,l)}u$ is in the image of $\rho$ 
for all $(\mu,k),(\nu,l)\in\lambda$ with $d^*(\mu)=d^*(\nu)$. 
Since the \Ca $A\otimes C(\T)$ is generated by those elements, 
the \shom $\rho$ is surjective. 
Therefore $\rho$ is an isomorphism between $\cT(F\times_{n,m}\T)$ 
and $A\otimes C(\T)$. 
This completes the proof of Lemma \ref{Lem:forappendix} 
as well as the one of Lemma \ref{Lem:appendix}. 

\begin{remark}\label{lastRemark}
With more efforts, 
we can show that 
the \Csa of $\cO(E\times_{n,m}\T)$ 
generated by $t^0(C(F^0\times\T))$ and $t^1(C(F^1\times\T))$ 
is isomorphic to 
\[
\bigoplus_{v\in S_1}\M_{k_v}\oplus
\bigoplus_{v\in S_2}\M_{k_v}\otimes C(\T) 
\]
where $S_1,S_2\subset F^0$ are defined by 
\begin{align*}
S_1&=\{v\in F^0\cap E^0_m\mid 
\emptyset\neq r^{-1}(v)\setminus F^1\subset m^{-1}(0)\}\\
S_2&=F^0\setminus \{v\in F^0\cap E^0_m\mid 
r^{-1}(v)\setminus F^1\subset m^{-1}(0)\}
\end{align*}
and $k_v$ is the number of the set $\{(\mu,k)\in \lambda\mid d^*(\mu)=v\}$
for $v\in S_1\cup S_2$. 
\end{remark}

\end{document}